\documentclass[11pt,a4paper]{art-num}
\usepackage{amsmath}
\usepackage{amssymb}
\usepackage{amsfonts}
\usepackage{hyperref}
\usepackage[dvips]{graphicx}

\newtheorem{theorem}{Theorem}

\newtheorem{definition}[theorem]{Definition}

\newtheorem{proposition}[theorem]{Proposition}

\newenvironment{proof}[1][Proof]{\noindent\textbf{#1.} }{\ \rule{0.5em}{0.5em}}
\renewcommand{\theequation}{\thesection.\arabic{equation}}
\input{tcilatex}

\def\Qlb#1{#1}

\def\Qcb#1{#1} 

\def\FRAME#1#2#3#4#5#6#7#8
{
 \begin{figure}[ptbh]
 \begin{center}
 \includegraphics[height=#3]{#7}
 \caption{#5}
 \label{#6}
 \end{center}
 \end{figure}
}

\reportno{AEI-2009-119}

\begin{document}

\title{Moduli spaces of $G_{2}$ manifolds}
\author{Sergey Grigorian \\
Max-Planck-Institut f\"{u}r Gravitationsphysik (Albert-Einstein-Institut)\\
Am M\"{u}hlenberg 1\\
D-14476 Golm\\
Germany}
\maketitle

\begin{abstract}
This paper is a review of current developments in the study of moduli spaces
of $G_{2}$ manifolds. $G_{2}$ manifolds are $7$-dimensional manifolds with
the exceptional holonomy group $G_{2}.$ Although they are odd-dimensional,
in many ways they can be considered as an analogue of Calabi-Yau manifolds
in $7$ dimensions. They play an important role in physics as natural
candidates for supersymmetric vacuum solutions of $M$-theory
compactifications. Despite the physical motivation, many of the results are
of purely mathematical interest. Here we cover the basics of $G_{2}$
manifolds, local deformation theory of $G_{2}$ structures and the local
geometry of the moduli spaces of $G_{2}$ structures.
\end{abstract}

\section{Introduction}

Ever since antiquity there has been a very close relationship between
physics and geometry. Originally, in \emph{Timaeus}, Plato related four of
the five Platonic solids - tetrahedron, hexahedron, octahedron, icosahedron
to the elements fire, earth, air and water, respectively, while the fifth
solid, the dodecahedron was the \emph{quintessence} of which the cosmos
itself is made. Later, Isaac Newton's Laws of Motion and Theory of
Gravitation gave a precise mathematical framework in which the motion of
objects can be calculated. However Albert Einstein's General Relativity made
it very explicit that the physics of spacetime is determined by its
geometry. More recently, this fundamental relationship has been taken to a
new level with the development of String and M-theory. Over the past 25
years, superstring theory has emerged as a successful candidate for the role
of a theory that would unify gravity with other interactions. It was later
discovered that all five superstring theories can be obtained as special
limits of a more general eleven-dimensional theory known as M-theory and
moreover, the low energy limit of which is the eleven-dimensional
supergravity \cite{Townsend:1995kk,Witten:1995ex}. The complete formulation
of M-theory is, however, not known yet.

One of the key features of String and M-theory is that these theories are
formulated in ten- and eleven-dimensional spacetimes, respectively. One of
the techniques to relate this to the visible four-dimensional world is to
assume that the remaining six or seven dimensions are curled up as a small,
compact, so-called \emph{internal space}. This is known as \emph{%
compactification}. Such a procedure also leads to a remarkable
interrelationship between physics and geometry, since the effective physical
content of the resulting four-dimensional theory is determined by the
geometry of the internal space. Usually the full multidimensional spacetime
is regarded as a direct product $M_{4}\times X$, where $M_{4}$ is a $4$%
-dimensional non-compact manifold with Lorentzian signature $\left(
-+++\right) $ and $X$ is a compact six or seven dimensional Riemannian
manifold. In general, the parameters that define the geometry of the
internal space give rise to massless scalar fields known as \emph{moduli},
and the properties of the moduli space are determined by the class of spaces
used in the compactification.

The properties of the internal space in String and M-theory
compactifications are governed by physical considerations. A key ingredient
of these theories is \emph{supersymmetry }\cite{WestSUSY}. Supersymmetry is
a physical symmetry between particles the spin of which differs by $\frac{1}{%
2}$ - that is, between integer spin \emph{bosons }and half-integer spin 
\emph{fermions}. Mathematically, bosons are represented as functions or
tensors and fermions as spinors. When looking for a supersymmetric vacuum
for which the metric is the only non-zero field, that is a Ricci-flat
solution that is invariant under supersymmetry transformations, it turns out
that a necessary requirement is the existence of covariantly constant, or
parallel, spinor. That is, there must exist a non-trivial spinor $\eta $ on
the Riemannian manifold $X$ that satisfies%
\begin{equation}
\nabla \eta =0  \label{constspineq}
\end{equation}%
where $\nabla $ is the relevant spinor covariant derivative \cite%
{BookSchwarzBecker}. This condition implies that $\eta $ is invariant under
parallel transport.

Properties of parallel transport on a Riemannian manifold are closely
related to the concept of \emph{holonomy}. Consider a vector $v$ at some
point $x$ on $X$. Using the natural Levi-Civita connection that comes from
the Riemannian metric, we can parallel transport $v$ along paths in $X$. In
particular, consider a closed contractible path $\gamma $ based at $x$. As
shown in Figure \ref{holofig}, if we parallel transport $v$ along $\gamma $,
then the new vector $v^{\prime }$ which we get will necessarily have the
same magnitude as the original vector $v$, but otherwise it does not have to
be the same. This gives the notion of \emph{holonomy group}. Below we give
the precise definition. 

\begin{definition}
Let $\left( X,g\right) $ be a Riemannian manifold of dimension $n$ with
metric $g$ and corresponding Levi-Civita connection $\nabla $, and fix point 
$x\in X$. Let $\gamma :\left[ 0,1\right] \longrightarrow X$ be a loop based
at $x$, that is, a piecewise-smooth path such that $\gamma \left( 0\right)
=\gamma \left( 1\right) =x.$ The parallel transport map $P_{\gamma
}:T_{x}X\longrightarrow T_{x}X$ is then an invertible linear map which lies
in $SO\left( n\right) $. Define the Riemannian \emph{holonomy group }$%
Hol_{x}\left( X,g\right) $ of $\nabla $ based at $x$ to be 
\begin{equation*}
Hol_{x}\left( X,g\right) =\left\{ P_{\gamma }:\gamma \ \text{is a loop based
at }x\right\} \subset O\left( n\right) 
\end{equation*}
\end{definition}

\FRAME{ftbpFU}{2.6085in}{2.298in}{0pt}{\Qcb{Parallel transport of a vector}}{%
\Qlb{holofig}}{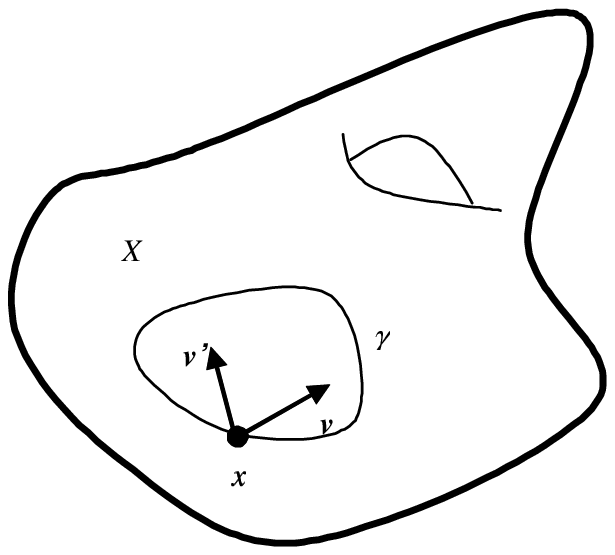}{\special{language "Scientific Word";type
"GRAPHIC";maintain-aspect-ratio TRUE;display "USEDEF";valid_file "F";width
2.6085in;height 2.298in;depth 0pt;original-width 2.4641in;original-height
2.1677in;cropleft "0";croptop "1";cropright "1";cropbottom "0";filename
'holonomy.eps';file-properties "XNPEU";}}

If the manifold $X$ is connected, then it is trivial to see that the
holonomy group is independent of the base point, and can hence be defined
for the whole manifold. Parallel transport is initially defined for vectors,
but can then be naturally extended to other objects like tensors and
spinors, with the holonomy group acting on these objects via relevant
representations.

Now going back to the covariantly constant spinor $\eta $, (\ref{constspineq}%
) implies that $\eta $ is invariant under the action of the holonomy group.
This shows that the spinor representation of $Hol\left( X,g\right) $ must
contain the trivial representation. For $Hol\left( X,g\right) =SO\left(
n\right) $, this is not possible since the spinor representation is
reducible, so $Hol\left( X,g\right) \subset SO\left( n\right) $. Hence the
condition (\ref{constspineq}) implies a reduced holonomy group. Thus,
Ricci-flat special holonomy manifolds occur very naturally in string and
M-theory.

As shown by Berger \cite{Berger1955}, the list of possible special holonomy
groups is very limited. In particular, if $X$ is simply-connected, and
neither locally a product nor symmetric, the only possibilities are given in
Figure \ref{hologroups}.%
\begin{figure}%
\begin{equation*}
\begin{tabular}{|l|l|l|}
\hline
\textbf{Geometry} & \textbf{Holonomy} & \textbf{Dimension} \\ \hline
K\"{a}hler & $U\left( k\right) $ & $2k$ \\ \hline
Calabi-Yau & $SU\left( k\right) $ & $2k$ \\ \hline
HyperK\"{a}hler & $Sp\left( k\right) $ & $4k$ \\ \hline
Exceptional & $G_{2}$ & $7$ \\ \hline
Exceptional & $Spin\left( 7\right) $ & $8$ \\ \hline
\end{tabular}%
\end{equation*}%
\caption{List of special holonomy groups}%
\label{hologroups}%
\end{figure}
In this list manifolds with holonomy $SU\left( k\right) ,Sp\left( k\right) ,$
$G_{2}$ and $Spin\left( 7\right) $ are Ricci-flat. Moreover, these groups
are subgroups of $SO\left( n\right) $ and are simply-connected. This implies
that manifolds with these holonomy groups always admit a spin structure (%
\cite[Proposition 3.6.2]{Joycebook}). These are also precisely the manifolds
that admit a parallel spinor. K\"{a}hler manifolds only admit parallel
projective spinors - a line subbundle of the spinor bundle. Thus, for a
Ricci-flat supersymmetric vacuum in a $10$-dimensional theory, $X$ has to be 
$6$-dimensional in order to reduce to $4$ dimensions, and hence necessarily
a Calabi-Yau manifold. Similarly, for an $11$-dimensional theory, $7$%
-dimensional manifolds with $G_{2}$ holonomy arise naturally. 

We have thus seen that even rather simple physical requirement restrict the
geometry of the manifold $X$ to rather special classes. In particular, the
study of Calabi-Yau manifolds has been crucial in the development of String
Theory, and in fact some very important discoveries in the theory of
Calabi-Yau manifolds have been made thanks to advances in the physics. One
such major discovery is Mirror Symmetry \cite%
{Strominger:1996it,MirrorSymBook}. This symmetry first appeared in String
Theory where evidence was found that conformal field theories (CFTs) related
to compactifications on a Calabi-Yau manifold with Hodge numbers $\left(
h_{1,1},h_{2,1}\right) $ are equivalent to CFTs on a Calabi-Yau manifold
with Hodge numbers $\left( h_{2,1},h_{1,1}\right) $. Mirror symmetry is
currently a powerful tool both for calculations in String Theory and in the
study of the Calabi-Yau manifolds and their moduli spaces.

In mathematical literature $G_{2}$ holonomy first appeared in Berger's list
of special holonomy groups in 1955 \cite{Berger1955}. In 1966 Bonan has
shown that manifolds with $G_{2}$ holonomy are Ricci-flat. It was known from
general theory that having a holonomy group $G$ is equivalent to having a
torsion-free $G$-structure. So it was natural to study $G_{2}$ structures on
manifolds to get a better understanding of $G_{2}$ holonomy. The different
classes of $G_{2}$ structures have been explored by Fern\'{a}ndez and Gray
in their 1982 paper \cite{FernandezGray}. In particular they have shown that
a torsion-free $G_{2}$ structure is equivalent to the $G_{2}$-invariant $3$%
-form $\varphi $ being closed and co-closed.

It was not known whether the group $G_{2}$ (or indeed $Spin\left( 7\right) $
for that matter) does actually appear as a non-symmetric holonomy group
until in 1987 Bryant \cite{Bryant-1987} proved the existence of metrics with 
$G_{2}$ and $Spin\left( 7\right) $ holonomy. In a later paper, Bryant and
Salamon \cite{BryantSalamon} constructed complete metrics with $G_{2}$
holonomy. However the first \emph{compact }examples of $G_{2}$ holonomy
manifolds have been constructed by Joyce in 1996 \cite{JoyceG2}. These
examples are based on quotients $T^{7}/\Gamma $ where $\Gamma $ is a finite
group. Such quotient spaces usually exhibit singularities, and Joyce has
shown that it is possible to resolve these singularities in such a way as to
get a smooth, compact manifold with $G_{2}$ holonomy. Since then, a number
of other types of constructions have been found, in particular the
construction by Kovalev \cite{Kovalev:2001zr} where a compact $G_{2}$
manifold is obtained by gluing together two non-compact asymptotically
cylindrical Riemannian manifolds with holonomy $SU\left( 3\right) $.

In the $G_{2}$ holonomy compactification approach to M-theory, the physical
content of the four-dimensional theory is given by the moduli of $G_{2}$
holonomy manifolds. Such a compactification of M-theory is in many ways
analogous to Calabi-Yau compactifications in String Theory, where much
progress has been made through the study of the Calabi-Yau moduli spaces. In
particular, as it was shown in \cite{Candelas:1990pi} and \cite%
{Strominger:1990pd}, the moduli space of complex structures and the
complexified moduli space of K\"{a}hler structures are both in fact, K\"{a}%
hler manifolds. Moreover, both have a \emph{special geometry}: that is, both
have a line bundle whose first Chern class coincides with the K\"{a}hler
class. However, until recently, the structure of the moduli space of $G_{2}$
holonomy manifolds has not been studied in that much detail. Generally, it
turns out that the study of $G_{2}$ manifolds is quite difficult. Unlike the
study of Calabi-Yau manifolds where the machinery of algebraic geometry has
been used with great success, in the case of $G_{2}$ manifolds there is no
analogue, so analytical rather than algebraic study is needed.

In this review, we aim to give an overview of what is currently known about $%
G_{2}$ moduli spaces and corresponding deformations of $G_{2}$ structures.
We first give an introduction to the properties of the group $G_{2}$ -
definitions and representations. Then we look at general properties of $%
G_{2} $ structures. Finally we move on to properties of $G_{2}$ moduli
spaces.

\section{The group $G_{2}$}

\setcounter{equation}{0}

\subsection{Automorphisms of octonions}

\label{sectg2oct}The group $G_{2}$ is the smallest of the $5$ exceptional
Lie groups, the others being $F_{4},E_{6},E_{7}$ and $E_{8}$. Surprisingly,
all of these Lie groups are related to the octonions, but $G_{2}$ is
especially close. So let us first give a few facts about the octonions. The
eight-dimensional algebra of octonions, denoted by $\mathbb{O}$, is the
largest possible normed division algebra. The others of course are the real
numbers $\mathbb{R}$, complex numbers $\mathbb{C}$ and the quaternions $%
\mathbb{H}$. Following Baez \cite{BaezOcto}, it turns out that division
algebras can be defined using the notion of \emph{triality}. Given three
real vector spaces $U,V,W$, then a triality is a non-degenerate trilinear
map 
\begin{equation*}
t:U\times V\times W\longrightarrow \mathbb{R}.
\end{equation*}%
Non-degenerate here means that for any fixed non-zero elements of $U$ and $V$%
, the induced functional on $W$ is non-zero. Hence, $t$ also defines a
bilinear map $m$ 
\begin{equation*}
m:U\times V\longrightarrow W^{\ast }.
\end{equation*}%
For each fixed element of $U$, this map defines an isomorphism between $V$
and $W^{\ast }$, and for each fixed element of $V$, an isomorphism between $%
U $ and \ $W^{\ast }$. Hence these three spaces are isomorphic to each, and
if we choose to identify non-zero elements $e_{1}\in U$, $e_{2}\in V$, and $%
e_{1}e_{2}\in W^{\ast }$, we can identify the spaces $U,V,W$ with each
other, and we can say that $m$ now defines multiplication on $U$ with
identity element $e=e_{1}=e_{2}=e_{1}e_{2}$. Note that in particular, the
existence of a non-degenerate trilinear map implies that the original vector
spaces $U$,$V$,$W$ are all of the same dimension.

Due to the non-degeneracy of the original triality, multiplication by a
fixed element is an isomorphism, so in fact, $U$ is a division algebra.
Assuming further that $U,V,W$ are inner product spaces, if the triality map
satisfies 
\begin{equation*}
\left\vert t\left( u,v,w\right) \right\vert \leq \left\Vert u\right\Vert
\left\Vert v\right\Vert \left\Vert w\right\Vert 
\end{equation*}%
and is such that for all $u,v$ there exists a non-zero $w$ such that the
bound is attained (and similarly for cyclic permutations for $u\,,v,w$) then
we get a normed division algebra. The converse is also true - any division
algebra defines a triality.

As discussed in detail by Baez \cite{BaezOcto}, on $\mathbb{R}^{n}$ it is
possible to construct bilinear maps $m_{n}$ involving the vector and spinor
representations of $Spin\left( n\right) $ 
\begin{subequations}%
\label{vecspinbm} 
\begin{subequations}
\begin{eqnarray}
m_{n} &:&V_{n}\times S_{n}^{\pm }\longrightarrow S_{n}^{\mp }\ \text{for }%
n=0,4\ \func{mod}8 \\
m_{n} &:&V_{n}\times S_{n}\longrightarrow S_{n}\mathbb{\ \ }\text{otherwise}
\end{eqnarray}%
\end{subequations}%
where $V_{n}$ is the vector representation of $SO\left( n\right) $, $%
S_{n}^{(\pm )}$ are the (left- and right-handed) spinor representations.

The spinor representations in (\ref{vecspinbm}) are self-dual, so in
principle, by dualizing the maps in (\ref{vecspinbm}), we could obtain
trilinear maps into $\mathbb{R}$. However, in order to obtain trialities,
these maps have to be non-degenerate, and hence the dimensions of the
relevant representations must agree. This happens only for $n=1,2,4,8,$ and
each of these trialities gives a normed division algebra of the
corresponding dimension: 
\end{subequations}
\begin{equation}
\begin{array}{cc}
t_{1}:V_{1}\times S_{1}\times S_{1}\longrightarrow \mathbb{R} & 
\Longrightarrow \mathbb{R} \\ 
t_{2}:V_{2}\times S_{2}\times S_{2}\longrightarrow \mathbb{R} & 
\Longrightarrow \mathbb{C} \\ 
t_{4}:V_{4}\times S_{4}^{+}\times S_{4}^{-}\longrightarrow \mathbb{R} & 
\Longrightarrow \mathbb{H} \\ 
t_{8}:V_{8}\times S_{8}^{+}\times S_{8}^{-}\longrightarrow \mathbb{R} & 
\Longrightarrow \mathbb{O}%
\end{array}%
\end{equation}%
This way, via the trialities we obtain all of the normed division algebras.

In general, suppose we have a triality $t:U_{1}\times U_{2}\times
U_{3}\longrightarrow \mathbb{R}$. Then to define a normed division algebra
from $t$, we fix two vectors in the two of the three spaces. Hence the
automorphism of the division algebra is the subgroup of the automorphism
group of the triality that fixes these two vectors. For $t_{8}$ the
automorphism group of the triality turns out to be $Spin\left( 8\right) $,
while $G_{2}$ is defined as the automorphism group of the corresponding
octonion algebra. Thus we have

\begin{definition}
The group $G_{2}$ is the automorphism group of the octonion algebra.
\end{definition}

Since $G_{2}$ is the automorphism group of octonions, it is the subgroup of $%
Spin\left( 8\right) $ (the automorphism group of the triality $t_{8}$) that
preserves unit vectors in $V_{8}$ and $S_{8}^{+}$. As explained by Baez in 
\cite{BaezOcto}, the subgroup of $Spin\left( 8\right) $ that fixes a unit
vector in $V_{8}$ is $Spin\left( 7\right) $. Moreover, if the representation 
$S_{8}^{+}$ is restricted to $Spin\left( 7\right) $, we get the spinor
representation $S_{7}$. Therefore, $G_{2}$ is the subgroup of $Spin\left(
7\right) $ that fixes a unit vector in $S_{7}$. In this representation, $%
Spin\left( 7\right) $ acts transitively on the unit sphere $S^{7}$, so we
have 
\begin{equation}
Spin\left( 7\right) /G_{2}=S^{7}.  \label{spin7g2}
\end{equation}%
Hence we have the following result.

\begin{proposition}
The group $G_{2}$ has dimension $14$.
\end{proposition}

\begin{proof}
From (\ref{spin7g2}), 
\begin{equation*}
\dim G_{2}=\dim \left( Spin\left( 7\right) \right) -\dim S^{7}=21-7=14.
\end{equation*}
\end{proof}

The automorphism group fixes the identity, so in fact $G_{2}$ acts
non-trivially on octonions that are orthogonal to the identity - the
imaginary octonions, denoted by $\func{Im}\left( \mathbb{O}\right) $ and
thus we get a natural $7$-dimensional representation of $G_{2}.$ A closer
look at this representation reveals another description of $G_{2}$. Using
octonion multiplication, we can define a cross product on $\func{Im}\left( 
\mathbb{O}\right) $ by 
\begin{equation}
a\times b=\func{Im}\left( ab\right) =\frac{1}{2}\left( ab-ba\right) .
\label{7dimcross}
\end{equation}%
But $G_{2}$ preserves octonion multiplication, hence any element of $G_{2}$
preserves the $7$-dimensional cross product. Alternatively, (\ref{7dimcross}%
) can be written as 
\begin{equation}
a\times b=ab+\left\langle a,b\right\rangle  \label{7dimcross1}
\end{equation}%
where $\left\langle ,\right\rangle $ is the octonionic inner product, in
general defined by 
\begin{equation*}
\left\langle a,b\right\rangle =\frac{1}{2}\left( a^{\ast }b+ba^{\ast
}\right) .
\end{equation*}%
Also, it can be shown that 
\begin{equation}
\left\langle a,b\right\rangle =-\frac{1}{6}\func{Tr}\left( a\times \left(
b\times \cdot \right) \right)  \label{7diminprodcross}
\end{equation}%
Therefore, from (\ref{7dimcross1}), multiplication of imaginary octonions
can be defined in terms of the cross product, hence any transformation
preserving the cross product preserves multiplication on $\func{Im}\left( 
\mathbb{O}\right) $, and is thus in $G_{2}.$ So, $G_{2}$ is precisely the
group that preserves the $7$-dimensional cross product.

Moreover, from the cross product we can form a \textquotedblleft scalar
triple product\textquotedblright\ on $\func{Im}\left( \mathbb{O}\right) $
given by 
\begin{equation}
\varphi _{0}\left( a,b,c\right) =\left\langle a,b\times c\right\rangle
=\left\langle a,bc\right\rangle .  \label{phivcp}
\end{equation}%
This defines $\varphi _{0}$ as an anti-symmetric trilinear functional - that
is, a $3$-form on $\mathbb{R}^{7}$. Equivalently, for a basis $e_{i}$ of $%
\func{Im}\left( \mathbb{O}\right) $, 
\begin{equation}
e_{i}\times e_{j}=\varphi _{0\ ij}^{\ k}e_{k}.  \label{phibrack}
\end{equation}%
So in this description, the components of $\varphi _{0}$ are essentially the
structure constants of the algebra of imaginary octonions.

A well-known way to encode the multiplication rules for the octonions is the 
\emph{Fano plane} \cite{BaezOcto}. It is shown in Figure \ref{fanoplane}. In
the diagram, the vertices $e_{1},...,e_{7}$ are the seven square roots of $%
-1 $. Multiplication follows along the six straight lines (sides of the
triangle and the altitudes) and along the central circle in the direction of
the arrows. So if $e_{i}$, $e_{j}$, $e_{k}$ are in this order on a straight
line, then $e_{i}e_{j}=e_{k}$ and $e_{j}e_{i}=-e_{k}$.

\FRAME{ftbpFU}{4.1154in}{3.6339in}{0pt}{\Qcb{Fano plane}}{\Qlb{fanoplane}}{%
fano}{\special{language "Scientific Word";type
"GRAPHIC";maintain-aspect-ratio TRUE;display "USEDEF";valid_file "F";width
4.1154in;height 3.6339in;depth 0pt;original-width 4.9199in;original-height
4.3429in;cropleft "0";croptop "1";cropright "1";cropbottom "0";filename
'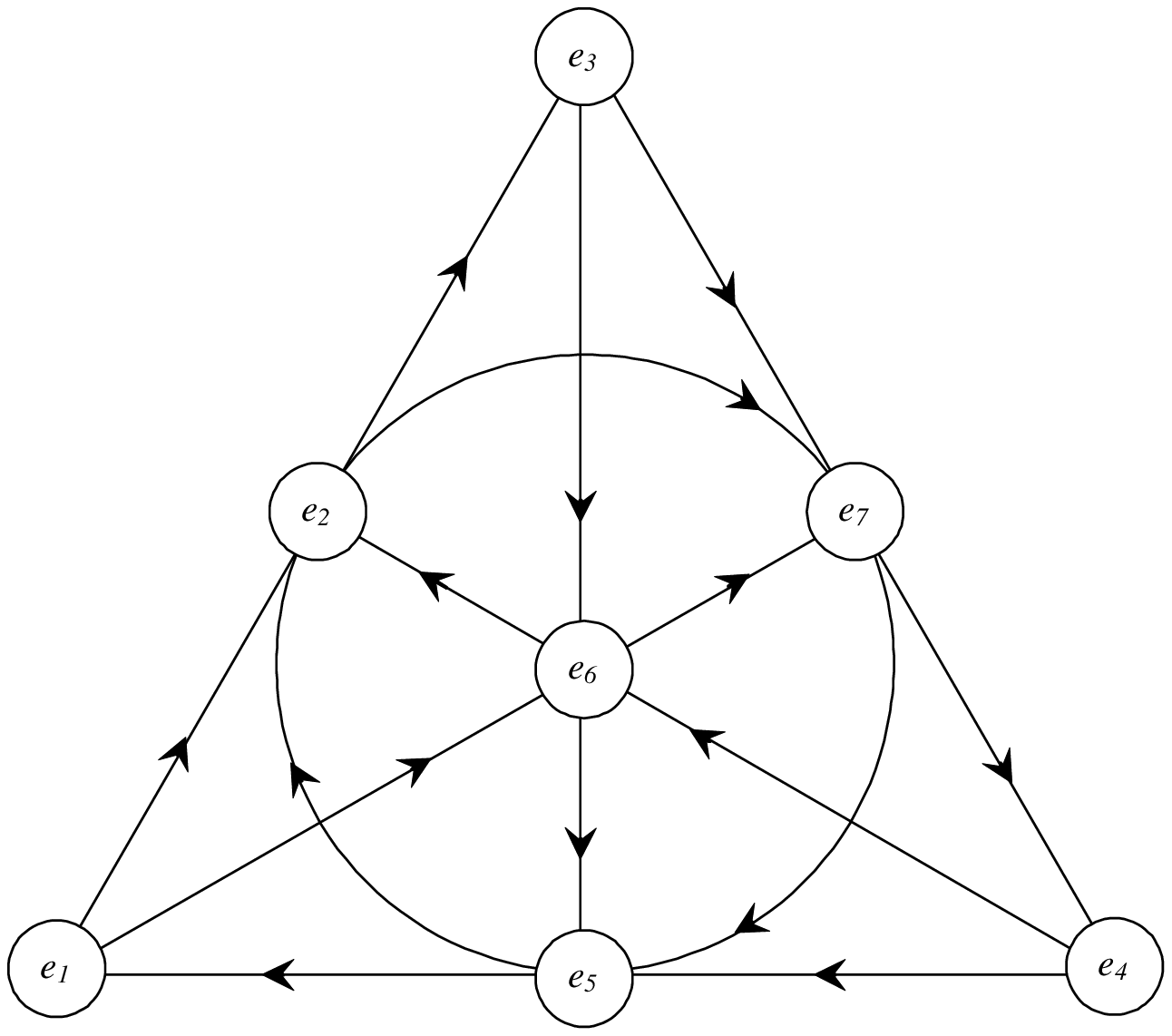';file-properties "XNPEU";}}

However, from (\ref{phibrack}) we see that $\varphi _{0}$ encodes precisely
the same information as the Fano plane. Suppose $x^{1},...,x^{7}$ are
coordinates on $\mathbb{R}^{7}$ and let $e^{ijk}=dx^{i}\wedge dx^{j}\wedge
dx^{k}$, then just reading off from the Fano plane, $\varphi _{0}$ can be
written as 
\begin{equation}
\varphi _{0}=e^{123}+e^{145}+e^{167}+e^{246}-e^{257}-e^{347}-e^{356}.
\label{phi0def}
\end{equation}%
Note that in order to keep the same convention for $\varphi _{0}$ as Joyce 
\cite{Joycebook}, in the Fano plane we have a different numbering for the
octonions compared to Baez \cite{BaezOcto}.

With this choice of coordinates, the inner product on $\func{Im}\left( 
\mathbb{O}\right) \cong \mathbb{R}^{7}$ is given by the standard Euclidean
metric 
\begin{equation}
g_{0}=\left( dx^{1}\right) ^{2}+...+\left( dx^{7}\right) ^{2}.  \label{g0def}
\end{equation}%
As seen from (\ref{7diminprodcross}), $G_{2}$ preserves the inner product on 
$\func{Im}\left( \mathbb{O}\right) $, so it clearly preserves $g_{0}$ and is
hence a subgroup of $SO\left( 7\right) $.

Since $\varphi _{0}$ defines the $7$-dimensional cross product, and $G_{2}$
is the symmetry group of this cross product, $G_{2}$ is the stabilizer of $%
\varphi _{0}$ in $GL\left( 7,\mathbb{R}\right) $. So we can state :

\begin{theorem}[Bryant, \protect\cite{Bryant-1987}]
\nocite{Harvey} The subgroup of $GL\left( 7,\mathbb{R}\right) $ that
preserves the $3$-form $\varphi _{0}$ is $G_{2}$.From the metric $g_{0}$ we
can define the Hodge star $\ast _{0}$ on $\mathbb{R}^{7}$, and using this,
the dual $4$-form $\psi _{0}=\ast _{0}\varphi _{0}$ which is given by 
\begin{equation}
\psi _{0}=e^{4567}+e^{2367}+e^{2345}+e^{1357}-e^{1346}-e^{1256}-e^{1247}.
\label{sphi0def}
\end{equation}
\end{theorem}

This is a key property of $G_{2}$ and as such this is often taken as the
definition of the group $G_{2}$, in particular in \cite{Joycebook}. As we
have seen, $G_{2}$ preserves both $\varphi _{0}$ and $g_{0}$, so it also
preserves $\psi _{0}$. In particular, $\varphi _{0}$ and $\psi _{0}$ give
alternate descriptions of the trivial $1$-dimensional representation of $%
G_{2}$.

It also turns out that $\psi _{0}$ is closely related to the \emph{%
associator }on $\func{Im}\left( \mathbb{O}\right) $. As the octonions are
non-associative, we can define a non-trivial associator map 
\begin{equation*}
\left[ \cdot ,\cdot ,\cdot \right] :\func{Im}\left( \mathbb{O}\right) \times 
\func{Im}\left( \mathbb{O}\right) \times \func{Im}\left( \mathbb{O}\right)
\longrightarrow \func{Im}\left( \mathbb{O}\right)
\end{equation*}%
given by 
\begin{equation}
\left[ a,b,c\right] =a\left( bc\right) -\left( ab\right) c.  \label{octassoc}
\end{equation}%
Just as $\varphi _{0}$ is defined as a dualization of the cross product
using the inner product to obtain the map 
\begin{equation*}
\varphi _{0}:\func{Im}\left( \mathbb{O}\right) \times \func{Im}\left( 
\mathbb{O}\right) \times \func{Im}\left( \mathbb{O}\right) \longrightarrow 
\mathbb{R}
\end{equation*}%
so it turns out that up to a constant multiple the map 
\begin{equation*}
\psi _{0}:\func{Im}\left( \mathbb{O}\right) \times \func{Im}\left( \mathbb{O}%
\right) \times \func{Im}\left( \mathbb{O}\right) \times \func{Im}\left( 
\mathbb{O}\right) \longrightarrow \mathbb{R}
\end{equation*}%
is a dualization of the associator, given by%
\begin{equation}
\psi _{0}\left( a,b,c,d\right) =\frac{1}{2}\left\langle \left[ a,b,c\right]
,d\right\rangle .  \label{psiassoc}
\end{equation}

It is possible to show that $\varphi _{0}$ and $\psi _{0}$ satisfy various
contraction identities. In particular, from \cite%
{bryant-2003,GrigorianYau1,karigiannis-2007}, we have

\begin{proposition}
The $3$-form $\varphi _{0}$ and the corresponding $4$-form $\psi _{0}$
satisfy the following identities: 
\begin{subequations}%
\label{contids} 
\begin{eqnarray}
\varphi _{0abc}\varphi _{0mn}^{\ \ \ \ c}
&=&g_{0am}g_{0bn}-g_{0an}g_{0bm}+\psi _{0abmn}  \label{phiphi1} \\
\varphi _{0abc}\psi _{0mnp}^{\ \ \ \ \ \ c} &=&3\left( g_{0a[m}\varphi
_{0np]b}-g_{0b[m}\varphi _{0np]a}\right)   \label{phipsi} \\
\psi _{0abcd}\psi _{0}^{mnpq} &=&24\delta _{a}^{[m}\delta _{b}^{n}\delta
_{c}^{p}\delta _{d}^{q]}+72\psi _{0[ab}^{\ \ \ \ [mn}\delta _{c}^{p}\delta
_{d]}^{q]}-16\varphi _{0[abc}\varphi _{0}^{\ [mnp}\delta _{d]}^{q]}
\label{psipsi0}
\end{eqnarray}%
\end{subequations}%
where $\left[ m\ n\ p\right] $ denotes antisymmetrization of indices and $%
\delta _{a}^{b}$ is the Kronecker delta, with $\delta _{b}^{a}=1$ if $a=b$
and $0$ otherwise.
\end{proposition}

The above identities can be of course further contracted - the details can
be found in \cite{GrigorianYau1,karigiannis-2007}. These identities and
their contractions are crucial whenever any calculations involving $\varphi
_{0}$ and $\psi _{0}$ have to be done. In particular, these are very useful
when studying $G_{2}$ manifolds.

\subsection{Representations of $G_{2}$}

As we will see in section \ref{g2structsec}, a crucial role in the study of $%
G_{2}$ structures is played by the representations of $G_{2}$. Since $G_{2}$
is a subgroup of $SO\left( 7\right) $, it has a fundamental vector
representation on $\mathbb{R}^{7}$ . In the study of $G_{2}$ manifolds, it
is very important to understand the representations of $G_{2}$ on $p$-forms.
So let us consider first the representations of $G_{2}$ on antisymmetric
tensors in $\mathbb{R}^{7}.$ For brevity let $V=\mathbb{R}^{7}.$ Following
Bryant \cite{bryant-2003}, we first look at the the Lie algebra $\mathfrak{so%
}\left( 7\right) $, which is the space of antisymmetric $7\times 7$ matrices
on $V$. For a vector $\omega \in V$, define the map

\begin{equation}
\rho _{\varphi }:V\longrightarrow \mathfrak{so}\left( 7\right) \ \text{given
by}\ \ \ \rho _{\varphi }\left( \omega \right) =\omega \lrcorner \varphi _{0}
\label{l27char}
\end{equation}%
which is clearly injective. Conversely, define the map 
\begin{equation}
\tau _{\varphi }:\mathfrak{so}\left( 7\right) \longrightarrow V\ \text{\
given by}\ \ \tau _{\varphi }\left( \alpha _{ab}\right) ^{c}=\frac{1}{6}%
\varphi _{0\ ab}^{\ c}\alpha ^{ab}.  \label{l214char}
\end{equation}%
From (\ref{contids}), we get that 
\begin{equation*}
\tau _{\varphi }\left( \rho _{\varphi }\left( \omega \right) \right) =\omega
,
\end{equation*}%
so that $\tau _{\varphi }$ is a partial inverse of $\rho _{\varphi }$. Thus
we get a decomposition 
\begin{equation}
\mathfrak{so}\left( 7\right) =\ker \tau _{\varphi }\oplus \rho _{\varphi
}\left( V\right)   \label{so7decom1}
\end{equation}%
where $\dim \rho _{\varphi }\left( V\right) =7$ and $\dim \ker \tau
_{\varphi }=14$. It turns out that $\ker \tau _{\varphi }$ is in fact a Lie
algebra with respect to the matrix commutator. This is the Lie algebra
bracket on $\mathfrak{so}\left( 7\right) $ and satisfies the Jacobi
identity. It is hence only necessary to show that for $\alpha ,\beta \in
\ker \tau _{\varphi },$ we have $\left[ \alpha ,\beta \right] \in \ker \tau
_{\varphi }$. This is an exercise in applying the contractions for $\varphi $%
.\ Thus we get a $14$-dimensional Lie subalgebra of $\mathfrak{so}\left(
7\right) $. However, this is precisely the Lie algebra $\mathfrak{g}_{2}$ 
\cite{karigiannis-2007}, that is%
\begin{equation}
\mathfrak{g}_{2}=\ker \tau _{\varphi }=\left\{ \alpha \in \mathfrak{so}%
\left( 7\right) :\varphi _{0abc}\alpha ^{bc}=0\right\} .  \label{liealgg2}
\end{equation}%
This further implies that we get the following decomposition of $\mathfrak{so%
}\left( 7\right) $:%
\begin{equation}
\mathfrak{so}\left( 7\right) =\mathfrak{g}_{2}\oplus \rho _{\varphi }\left(
V\right) .  \label{so7decom}
\end{equation}%
The group $G_{2}$ acts via the adjoint representation on the $14$%
-dimensional vector space $\mathfrak{g}_{2}$ and via the fundamental vector
representation on the $7$-dimensional space $\rho _{\varphi }\left( V\right) 
$. This is a $G_{2}\,$-invariant irreducible decomposition of $\mathfrak{so}%
\left( 7\right) $ into the representations $\mathbf{7}$ and $\mathbf{14}$.
Hence we get the following result:

\begin{theorem}[Bryant, \protect\cite{Bryant-1987}]
The space $\Lambda ^{2}$ of $2$-forms on $V$ decomposes as 
\begin{equation}
\Lambda ^{2}=\Lambda _{7}^{2}\oplus \Lambda _{14}^{2}.  \label{2formdecom}
\end{equation}%
with the components $\Lambda _{7}^{2}$ and $\Lambda _{14}^{2}$ given by: 
\begin{subequations}%
\label{lambda2decom} 
\begin{eqnarray}
\Lambda _{7}^{2} &=&\left\{ \omega \lrcorner \varphi :\omega \ \text{a vector%
}\right\}  \label{om27} \\
\Lambda _{14}^{2} &=&\left\{ \alpha =\frac{1}{2}\alpha _{ab}e^{a}\wedge
e^{b}:\left( \alpha _{ab}\right) \in \mathfrak{g}_{2}\right\}  \label{om214}
\end{eqnarray}%
\end{subequations}%
\end{theorem}

An alternative, but fully equivalent, description of $\Lambda _{7}^{2}$ and $%
\Lambda _{14}^{2}$ presents them as eigenspaces of the operator 
\begin{equation}
T_{\psi }:\Lambda ^{2}\longrightarrow \Lambda ^{2}\ \ \text{given by }%
T_{\psi }\left( \alpha _{ab}\right) =\psi _{0abcd}\alpha ^{cd}
\label{Tpsiop}
\end{equation}%
With this description, we have \cite{karigiannis-2007}: 
\begin{subequations}%
\label{tpsil} 
\begin{eqnarray}
\Lambda _{7}^{2} &=&\left\{ \alpha \in \Lambda ^{2}:T_{\psi }\alpha =4\alpha
\right\}  \label{om27a} \\
\Lambda _{14}^{2} &=&\left\{ \alpha \in \Lambda ^{2}:T_{\psi }\alpha
=-2\alpha \right\} \text{.}  \label{om214a}
\end{eqnarray}%
\end{subequations}%
Correspondingly, the description of the $\mathbf{7}$ and $\mathbf{14}$
pieces of $\Lambda ^{5}$ is obtained from (\ref{om27}) and (\ref{om214}) via
Hodge duality.

Let us now look at $3$-forms in more detail. Consider $\mathrm{Sym}%
^{2}\left( V^{\ast }\right) $ - the space of symmetric $2$-tensors on $V$,
and define a map 
\begin{equation}
\mathrm{i}_{\varphi }:\mathrm{Sym}^{2}\left( V^{\ast }\right)
\longrightarrow \Lambda ^{3}\ \ \text{given by }\mathrm{i}_{\varphi }\left(
h\right) _{abc}=h_{[a}^{d}\varphi _{0bc]d}.  \label{iphigen}
\end{equation}%
We can decompose $\mathrm{Sym}^{2}\left( V^{\ast }\right) =\mathbb{R}%
g_{0}\oplus \mathrm{Sym}_{0}^{2}\left( V^{\ast }\right) $ where $\mathbb{R}%
g_{0}$ is the set of symmetric tensors proportional to the metric $g_{0}$
and $\mathrm{Sym}_{0}^{2}\left( V^{\ast }\right) $ is the set of traceless
symmetric tensors. This is a $G_{2}$-invariant irreducible decomposition of $%
\mathrm{Sym}^{2}\left( V^{\ast }\right) $ into $1$-dimensional and $27$%
-dimensional representations. We clearly have 
\begin{equation*}
\mathrm{i}_{\varphi }\left( g_{0}\right) _{abc}=\varphi _{0abc},
\end{equation*}%
so the map $\mathrm{i}_{\varphi }$ is also $G_{2}$-invariant and is
injective on each summand of this decomposition. Looking at the first
summand, we get that $\mathrm{i}_{\varphi }\left( \mathbb{R}g_{0}\right)
=\Lambda _{1}^{3}$ - the one-dimensional singlet representation of $G_{2}.$
Now look at the second summand and consider $\mathrm{i}_{\varphi }\left( 
\mathrm{Sym}_{0}^{2}\left( V^{\ast }\right) \right) $. This is $27$%
-dimensional and irreducible, so it gives a $27$-dimensional representation
of $G_{2}$ on $3$-forms: 
\begin{equation*}
\mathrm{i}_{\varphi }\left( \mathrm{Sym}_{0}^{2}\left( V^{\ast }\right)
\right) =\Lambda _{27}^{3}\left( V^{\ast }\right) .
\end{equation*}%
Now, $\Lambda ^{3}$ is $35$-dimensional, and we have accounted for $1+27=28$
dimensions. Thus we still have $7$ dimensions left unaccounted for in $%
\Lambda ^{3}$. So let us extend the map $\mathrm{i}_{\varphi }$ to $\Lambda
^{2}$ - the antisymmetric $2$-tensors on $\mathbb{R}^{7}$. Suppose $\beta
\in \Lambda _{7}^{2}$. Then $\beta =\omega \lrcorner \varphi _{0}$, for some
vector $\omega \in V$ so 
\begin{equation}
\mathrm{i}_{\varphi }\left( \beta \right) _{abc}=\varphi _{0\ [a\left\vert
e\right\vert }^{d}\varphi _{0bc]d}^{\ }\omega ^{e}=\psi _{0abcd}\omega ^{d}
\label{iphil7}
\end{equation}%
where we have used (\ref{contids}). This defines a $G_{2}$-invariant map
from $V$ to $\Lambda ^{3}$ and hence gives $\Lambda _{7}^{3}$.

So overall we thus have a decomposition of $3$-forms into irreducible
representations of $G_{2}$:

\begin{theorem}[Bryant, \protect\cite{bryant-2003}]
The space $\Lambda ^{3}$ of $3$-forms on $V$ decomposes as%
\begin{equation}
\Lambda ^{3}=\Lambda _{1}^{3}\oplus \Lambda _{7}^{3}\oplus \Lambda _{27}^{3}
\label{lamb3dec}
\end{equation}%
where 
\begin{subequations}%
\label{lamb3decom} 
\begin{eqnarray}
\Lambda _{1}^{3} &=&\left\{ \chi \in \Lambda ^{3}:\chi _{abc}=f\varphi
_{0abc}\ \text{for scalar }f\right\}  \label{lamb31} \\
\Lambda _{7}^{3} &=&\left\{ \omega \lrcorner \psi _{0}:\omega \ \text{a
vector}\right\} .  \label{lamb37} \\
\Lambda _{27}^{3} &=&\left\{ \chi \in \Lambda ^{3}:\chi
_{abc}=h_{[a}^{d}\varphi _{0bc]d}\text{ for }h_{ab}~\text{traceless,
symmetric}\right\} .  \label{lamb327}
\end{eqnarray}%
\end{subequations}%
\end{theorem}

From the identities for contraction of $\varphi _{0}$ and $\psi _{0}$, it is
possible to see that an equivalent description of $\Lambda _{27}^{3}$ is 
\begin{equation*}
\Lambda _{27}^{3}=\left\{ \chi \in \Lambda ^{3}:\chi \wedge \varphi _{0}=0\ 
\text{and }\chi \wedge \psi _{0}=0\right\} .
\end{equation*}%
A similar decomposition of $4$-forms is again obtained via Hodge duality.

Suppose we have $\chi \in \Lambda ^{3}$, then define $\pi _{1}$, $\pi _{7}$
and $\pi _{27}$ to be projections of $\chi $ onto $\Lambda _{1}^{3}$, $%
\Lambda _{7}^{3}$ and $\Lambda _{27}^{3}$, respectively. Using contraction
identities for $\varphi $ and $\psi $, we get the following relations \cite%
{GrigorianYau1}:

\begin{proposition}
Given a $3$-form $\chi \in \Lambda ^{3}$, the projections of $\chi $ onto
the components (\ref{lamb3dec}) of $\Lambda ^{3}$ are given by: 
\begin{subequations}%
\label{3formproj} 
\begin{eqnarray}
\pi _{1}\left( \chi \right) &=&a\varphi _{0}\ \text{where }a=\frac{1}{42}%
\left( \chi _{abc}\varphi _{0}^{abc}\right) =\frac{1}{7}\,\left\langle \chi
,\varphi _{0}\right\rangle \ \text{with}\ \left\vert \pi _{1}\left( \chi
\right) \right\vert ^{2}=7a^{2} \\
\pi _{7}\left( \chi \right) &=&\omega \lrcorner \psi _{0}\ \text{where }%
\omega ^{a}=-\frac{1}{24}\chi _{mnp}\psi _{0}^{mnpa}\ \ \text{with \ }%
\left\vert \pi _{7}\left( \chi \right) \right\vert ^{2}=4\left\vert \omega
\right\vert ^{2} \\
\pi _{27}\left( \chi \right) &=&\mathrm{i}_{\varphi }\left( h\right) \ \text{%
where }h_{ab}=\frac{3}{4}\chi _{mn\{a}\varphi _{0b\}}^{\ \ mn}\ \text{with }%
\left\vert \pi _{27}\left( \chi \right) \right\vert ^{2}=\frac{2}{9}%
\left\vert h\right\vert ^{2}.
\end{eqnarray}%
\end{subequations}%
Here $\{a$ $b\}$ denotes the traceless symmetric part.
\end{proposition}

Note that similar projections can be defined for $4$-forms as well.

\section{$G_{2}$ structures}

\setcounter{equation}{0}\label{g2structsec}

\subsection{Definition}

As we shall see, the notion of holonomy is closely related to $\emph{G}$-%
\emph{structures }on manifolds. Let us give the necessary definitions

\begin{definition}
Let $X$ be a manifold of dimension $n$. Suppose $TX$ is the tangent bundle
over $X$. Define the manifold $F$ by%
\begin{equation*}
F=\left\{ \left( x,e_{1},...,e_{n}\right) :x\in X\ \ \text{and }\left(
e_{1},...,e_{n}\right) \ \text{is a basis for }T_{x}X\right\} 
\end{equation*}%
This then has a projection $\pi :\left( x,e_{1},...,e_{n}\right) \longmapsto
x$ onto $X$ and a natural left action by $GL\left( n,\mathbb{R}\right) $ on
the fibres. $F$ is thus a principal bundle over $X$ with fibre $GL\left( n,%
\mathbb{R}\right) ,$ called the \emph{frame bundle }of $X.$
\end{definition}

\begin{definition}
Let $X$ be a manifold of dimension $n$. Let $G$ be a Lie subgroup of $%
GL\left( n,\mathbb{R}\right) $. Then a $\emph{G}$-\emph{structure }on $X$ is
a principal subbundle $P$ of $F$ with fibre $G$.
\end{definition}

The framework of $G$-structures is very powerful, and a number of
geometrical structures can be reformulated in this language. In particular,
a Riemannian metric on a manifold is equivalent to an $O\left( n\right) $
structure. We are in particular interested in \emph{torsion-free }$G$%
-structures. A $G$-structure is torsion-free if and only if there exists a 
\emph{compatible }torsion-free connection on $TM$. A connection $\nabla $ on 
$TM$ is equivalent to a connection $D$ on the frame bundle $F$, and we say $%
\nabla $ is compatible with the $G$-structure $P$ if $D$ reduces to a
connection on $P$. For example, given a Riemannian metric, a unique
torsion-free Levi-Civita connection can always be defined, hence all $%
O\left( n\right) $ structures are torsion-free. On a complex manifold with
complex dimension, an integrable complex structure is equivalent to a
torsion-free $GL\left( m,\mathbb{C}\right) $ structure. A K\"{a}hler
structure is then equivalent to a torsion-free $U\left( m\right) $%
-structure. From \cite{Joycebook} we have a key result that relates
torsion-free structures and holonomy:

\begin{proposition}
\label{propholstruct}Let $\left( X,g\right) $ be a Riemannian manifold of
dimension $n$, with $O\left( n\right) $-structure $P$ corresponding to $g$.
Let $G$ be a Lie subgroup of $O\left( n\right) $. Then $Hol\left( g\right)
\subseteq G$ if and only if $X$ admits a torsion-free $G$-structure $Q$ that
is a subbundle of $P$.
\end{proposition}

As Proposition \ref{propholstruct} shows, the study of Riemannian holonomy
is equivalent to studying torsion-free $G$-structures. Hence in order to
study $G_{2}$ holonomy manifolds we will first consider $G_{2}$ structures.

Now suppose $X$ is a smooth, oriented $7$-dimensional manifold. Following
Joyce \cite{Joycebook}, define a $3$-form $\varphi $ to be \emph{positive }%
if locally we can choose a frame such that $\varphi $ is written in the form
(\ref{phi0def}) - that is for every $p\in X$ there is an oriented
isomorphism $q_{p}$ between $T_{p}X$ and $\mathbb{R}^{7}$ such that $\left.
\varphi \right\vert _{p}=\varphi _{0}$. For each $p\in X$ define $\mathcal{P}%
_{p}^{3}X$ to be set of such $3$-forms. To each positive $\varphi $ we can
associate a metric $g$ and a Hodge dual $\ast \varphi $ which are identified
with $g_{0}$ and $\psi _{0}$ under the $q_{p}$ and the associated metric is
written (\ref{g0def}).

Since $\varphi _{0}$ is preserved by $G_{2}$ and $GL\left( 7,\mathbb{R}%
\right) _{+}$ acts transitively on $\mathcal{P}_{p}^{3}X$ it follows that 
\begin{equation*}
\mathcal{P}_{p}^{3}X\cong GL\left( 7,\mathbb{R}\right) _{+}/G_{2}
\end{equation*}%
and hence $\dim \mathcal{P}_{p}^{3}X=\dim GL\left( 7,\mathbb{R}\right)
_{+}-\dim G_{2}=49-14=35$. This is equal to the dimension of $\Lambda
^{3}T_{p}^{\ast }X$, hence $\mathcal{P}_{p}^{3}X$ is an open subset of $%
\Lambda ^{3}T_{p}^{\ast }X$. Moreover if we consider the bundle $\mathcal{P}%
^{3}X$ over $X$ with fibre $\mathcal{P}_{p}^{3}X$, it will be an open
subbundle of $\Lambda ^{3}T^{\ast }X$.

Given a positive $3$-form $\varphi $ on $X$, consider at each point $p$ the
set $Q_{p}$ of isomorphisms $q_{p}$ between $T_{p}X$ and $\mathbb{R}^{7}$
such that $\left. \varphi \right\vert _{p}=\varphi _{0}$. It is then easy to
see that $Q_{p}$ $\cong G_{2}$ and that the bundle $Q$ over $X$ with fibre $%
Q_{p}$ is in fact a principal subbundle of the frame bundle $F$. So in fact, 
$Q$ is a $G_{2}$ structure. The converse is also true - given an oriented $%
G_{2}$ structure $Q$, we can uniquely define a positive $3$-form $\varphi $
and associated metric $g$ and $4$-form $\psi $ that correspond to $\varphi
_{0}$,$g_{0}$ and $\psi _{0}$ respectively. We thus have a key result:

\begin{theorem}[Joyce, \protect\cite{Joycebook}]
Let $X$ be an oriented $7$-dimensional manifold. There exists a $1-1$
correspondence between positive $3$-forms on $X$ and oriented $G_{2}$%
-structures $Q$ on $X$. Moreover, to each positive $3$-form $\varphi $ we
can associate a Riemannian metric $g$ and a corresponding $4$-form $\ast
_{\varphi }\varphi =\psi $ such for each $p\in X$, under the isomorphism $%
q_{p}:T_{p}X\longrightarrow \mathbb{R}^{7}$, these quantities are identified
with $\varphi _{0}$,$g_{0}$ and $\psi _{0}$ respectively.
\end{theorem}

So given a positive $3$-form $\varphi $ on $X$, it is possible to define a
metric $g$ associated to $\varphi $. This metric then defines the Hodge
star, which we denote by $\ast _{\varphi }$ to emphasize the dependence on $%
\varphi $. Given the Hodge star, we can in turn define the $4$-form $\psi
=\ast _{\varphi }\varphi $. Thus in fact both the metric $g$ and the $4$%
-form $\psi $ are functions of $\varphi $. By definition, at point $p\in X$
there is an isomorphism that identifies $\varphi $ with $\varphi _{0}$, $%
\psi $ with $\psi _{0}$ and $g$ with $g_{0}$. Therefore, properties of $%
\varphi _{0}$ and $\psi _{0}$ such as the contraction identities (\ref%
{contids}) that we encountered in Section \ref{sectg2oct} also hold for the
differential forms $\varphi $ and $\psi $.

In general, any $G$-structure on a manifold $X$ induces a splitting of
bundles of $p$-forms into subbundles corresponding to irreducible
representations of $G$. The same is of course true for $G_{2}\,$structures.
The decomposition of $p$-forms on $\mathbb{R}^{7}$ carries over to any
manifold with a $G_{2}$ structure, so from the previous section we have the
following decomposition of the spaces of $p$-forms $\Lambda ^{p}$: 
\begin{subequations}%
\label{formdecompose} 
\begin{eqnarray}
\Lambda ^{1} &=&\Lambda _{7}^{1}  \label{l1decom} \\
\Lambda ^{2} &=&\Lambda _{7}^{2}\oplus \Lambda _{14}^{2}  \label{l2decom} \\
\Lambda ^{3} &=&\Lambda _{1}^{3}\oplus \Lambda _{7}^{3}\oplus \Lambda
_{27}^{3}  \label{l3decom} \\
\Lambda ^{4} &=&\Lambda _{1}^{4}\oplus \Lambda _{7}^{4}\oplus \Lambda
_{27}^{4}  \label{l4decom} \\
\Lambda ^{5} &=&\Lambda _{7}^{5}\oplus \Lambda _{14}^{5}  \label{l5decom} \\
\Lambda ^{6} &=&\Lambda _{7}^{6}  \label{l6decom}
\end{eqnarray}%
\end{subequations}%

Here each $\Lambda _{k}^{p}$ corresponds to the $k\,$-dimensional
irreducible representation of $G_{2}$. Moreover, for each $k$ and $p$, $%
\Lambda _{k}^{p}$ and $\Lambda _{k}^{7-p}$ are isomorphic to each other via
Hodge duality, and also $\Lambda _{7}^{p}$ are isomorphic to each other for $%
n=1,2,...,6$.

Define the standard inner product on $\Lambda ^{p},$ so that for $p\,$-forms 
$\alpha $ and $\beta $, 
\begin{equation}
\left\langle \alpha ,\beta \right\rangle =\frac{1}{p!}\alpha
_{a_{1}...a_{p}}\beta ^{a_{1}...a_{p}}\text{.}  \label{forminp}
\end{equation}%
This is related to the Hodge star, since 
\begin{equation}
\alpha \wedge \ast \beta =\left\langle \alpha ,\beta \right\rangle \mathrm{%
vol}  \label{hodgedef}
\end{equation}%
where $\mathrm{vol}$ is the invariant volume form given locally by 
\begin{equation}
\mathrm{vol}=\sqrt{\det g}dx^{1}\wedge ...\wedge dx^{7}.  \label{voldef}
\end{equation}%
Then the decompositions (\ref{formdecompose}) are orthogonal with respect to
(\ref{forminp}). Note that $\left\langle \varphi ,\varphi \right\rangle =7$,
so in fact we have 
\begin{equation}
V=\frac{1}{7}\int \varphi \wedge \ast \varphi  \label{phiwpsi}
\end{equation}%
where $V$ is the volume of the manifold $X$.

We know that the metric $g$ is defined by the $3$-form $\varphi $ and we can
use some of the results from Section \ref{sectg2oct} to find a direct
relationship between the two quantities.

\begin{proposition}
Given a positive $3$-form $\varphi $ on a $7$-manifold $X$, the associated
metric $g$ is given by 
\begin{equation}
g_{ab}=\left( \det s\right) ^{-\frac{1}{9}}s_{ab}.  \label{metricdefdirect}
\end{equation}%
with 
\begin{equation}
s_{ab}=\frac{1}{144}\varphi _{amn}\varphi _{bpq}\varphi _{rst}\hat{%
\varepsilon}^{mnpqrst}  \label{sabdef}
\end{equation}%
where $\hat{\varepsilon}^{mnpqrst}$ is the alternating symbol with $\hat{%
\varepsilon}^{12...7}=+1$. Alternatively, for $u$,$v$ vector fields on $X$, 
\begin{equation}
\left\langle u,v\right\rangle \mathrm{vol}=\frac{1}{6}\left( u\lrcorner
\varphi \right) \wedge \left( v\lrcorner \varphi \right) \wedge \varphi
\label{metricdef}
\end{equation}%
where $\lrcorner $ denotes interior multiplication: $\left( u\lrcorner
\varphi \right) _{bc}=u^{a}\varphi _{abc}$.
\end{proposition}

\begin{proof}
Consider the quantity $P_{ab}$ given by 
\begin{equation*}
P_{ab}=\varphi _{amn}\varphi _{bpq}\psi ^{mnpq}
\end{equation*}%
Using identities (\ref{contids}) to contract $\varphi $ and $\psi $, this
gives 
\begin{equation*}
P_{ab}=24g_{ab}.
\end{equation*}%
Expanding $\psi ^{mnpq}$ in terms of $\varphi $ and the Levi-Civita tensor
we get 
\begin{equation*}
P_{ab}=\frac{1}{6}\varphi _{amn}\varphi _{bpq}\varphi _{rst}\varepsilon
^{mnpqrst}\text{.}
\end{equation*}%
If we write $\hat{\varepsilon}^{mnpqrst}$ for the alternating symbol with $%
\hat{\varepsilon}^{12...7}=+1$, then we get 
\begin{equation}
g_{ab}\sqrt{\det g}=\frac{1}{144}\varphi _{amn}\varphi _{bpq}\varphi _{rst}%
\hat{\varepsilon}^{mnpqrst}.  \label{metriccomp1}
\end{equation}%
Alternatively, let $u$ and $v$ be vector fields on $X$. Then%
\begin{equation*}
\left\langle u,v\right\rangle \sqrt{\det g}=\frac{1}{144}\left( u^{a}\varphi
_{amn}\right) \left( v^{a}\varphi _{bpq}\right) \varphi _{rst}\hat{%
\varepsilon}^{mnpqrst}.
\end{equation*}%
Hence we get (\ref{metricdef}). Now define%
\begin{equation*}
s_{ab}=\frac{1}{144}\varphi _{amn}\varphi _{bpq}\varphi _{rst}\hat{%
\varepsilon}^{mnpqrst}
\end{equation*}%
so that then, after taking the determinant of (\ref{metriccomp1}) we get (%
\ref{metricdefdirect}).
\end{proof}

Thus we see that even though given the $3$-form $\varphi $ we can define the
metric $g$, this relationship is rather complicated and non-linear. In
particular, this also shows that $\psi =\ast _{\varphi }\varphi $ depends on 
$\varphi $ in an even more non-trivial fashion, since the Hodge star depends
itself on the metric.

Here we need to say a few words about the notation used for the $G_{2}$ $3$%
-form $\varphi $ and the associated $4$-form $\psi $. The notation that we
use here is due to Karigiannis - where the Hodge dual of $\varphi $ is
denoted by $\psi $ and was first introduced in \cite{karigiannis-2005-57}.
In Figure \ref{notation} we summarize the different notations used by other
authors: 
\begin{figure}%

\begin{equation*}
\begin{tabular}{|l|l|l|l|}
\hline
\textbf{Authors} & $3$\textbf{-form} & \textbf{Dual }$4$\textbf{-form} & 
\textbf{References} \\ \hline
\begin{tabular}{l}
Beasley and Witten \\ 
Gukov, Yau and Zaslow%
\end{tabular}
& $\Phi $ & $\ast \Phi $ & \cite{WittenBeasley,Gukov:2002jv} \\ \hline
Bryant & $\phi =\frac{1}{6}\varepsilon _{ijk}e^{ijk}$ & $\ast _{\phi }\phi =%
\frac{1}{24}\varepsilon _{ijkl}e^{ijkl}$ & \cite{Bryant-1987,bryant-2003} \\ 
\hline
Hitchin; Lee and Leung & $\Omega $ & $\Theta =\ast \Omega $ & \cite%
{Hitchin:2000jd,Lee:2002fa} \\ \hline
Joyce & $\varphi $ & $\Theta \left( \varphi \right) =\ast \varphi $ & \cite%
{JoyceG2,Joycebook} \\ \hline
\begin{tabular}{l}
Karigiannis; Karigiannis and Leung \\ 
Grigorian and Yau%
\end{tabular}
& $\varphi $ & $\psi =\ast _{\varphi }\varphi $ & \cite%
{GrigorianYau1,karigiannis-2005-57,karigiannis-2007,karigiannis-2007a} \\ 
\hline
\end{tabular}%
\end{equation*}%
\caption{Notation that is used by different authors}%
\label{notation}%
\end{figure}%
where $e^{ijk}=e^{i}\wedge e^{j}\wedge e^{k}$ and $e^{ijkl}=e^{i}\wedge
e^{j}\wedge e^{k}\wedge e^{l}$ for basis covectors $e^{i}$.

\subsection{Torsion-free structures}

The definition of a $G_{2}$ structure only defines the algebraic properties
of $\varphi $, and in general does not address the analytical properties of $%
\varphi $. Using the associated metric $g$ we can define the Levi-Civita
connection $\nabla $ on $X$. Then it is natural to ask what are the
properties of $\nabla \varphi $. This quantity is known as the \emph{torsion}
of the $G_{2}$ structure. Originally the torsion of $G_{2}$ structures was
studied by Fern\'{a}ndez and Gray \cite{FernandezGray}, and their analysis
revealed that there are in fact a total of 16 torsion classes of $G_{2}$
structures. Later on, Karigiannis reproduced their results using simple
computational arguments \cite{karigiannis-2007}.

Following \cite{karigiannis-2007}, consider the $3$-form $\nabla _{X}\varphi 
$ for some vector field $X$. We know that $3$-forms split as $\Lambda
_{1}^{3}\oplus \Lambda _{7}^{3}\oplus \Lambda _{27}^{3}$, so consider the
projections $\pi _{1}$,$\pi _{7}$ and $\pi _{27}$ of $\nabla _{X}\varphi $
onto these components. Using (\ref{3formproj}), we have 
\begin{equation*}
\pi _{1}\left( \nabla _{X}\varphi \right) =a\varphi
\end{equation*}%
where 
\begin{eqnarray*}
a &=&X^{a}\left( \nabla _{a}\varphi _{bcd}\right) \varphi ^{bcd}=X^{a}\nabla
_{a}\left( \varphi _{bcd}\varphi ^{bcd}\right) -\varphi _{bcd}X^{a}\nabla
_{a}\varphi ^{bcd} \\
&=&-X^{a}\left( \nabla _{a}\varphi _{bcd}\right) \varphi ^{bcd} \\
&=&0.
\end{eqnarray*}%
Hence we see that the $\Lambda _{1}^{3}$ component vanishes. Similarly, for $%
\Lambda _{27}^{3}$ we have 
\begin{equation*}
\pi _{27}\left( \nabla _{X}\varphi \right) =\mathrm{i}_{\varphi }\left(
h\right) \ 
\end{equation*}%
where 
\begin{eqnarray*}
\text{ }h_{ab} &=&\frac{3}{4}\left( X^{c}\nabla _{c}\varphi _{mn\{a}\right)
\varphi _{b\}}^{\ \ mn}=\frac{3}{4}X^{c}\nabla _{c}\left( \varphi
_{mn\{a}\varphi _{b\}}^{\ \ \ mn}\right) -\frac{3}{4}\varphi
_{mn\{a}X^{c}\nabla _{c}\varphi _{b\}}^{\ \ \ mn} \\
&=&-\frac{3}{4}\left( X^{c}\nabla _{c}\varphi _{mn\{a}\right) \varphi
_{b\}}^{\ \ mn} \\
&=&0.
\end{eqnarray*}%
Here we have used the fact that $\varphi _{mna}\varphi _{b}^{\ \ \
mn}=6g_{ab}$, the traceless part of which vanishes. Therefore, the $\Lambda
_{27}^{3}$ part of $\nabla _{X}\varphi $ also vanishes. Now consider the $%
\Lambda _{7}^{3}$ component. In this case,%
\begin{equation*}
\pi _{7}\left( \nabla _{X}\varphi \right) =\omega \lrcorner \psi
\end{equation*}%
$\ $where 
\begin{equation*}
\omega ^{a}=-\frac{1}{24}X^{c}\left( \nabla _{c}\varphi _{mnp}\right) \psi
^{mnpa}=\frac{1}{24}X^{a}\left( \nabla _{a}\psi ^{bcde}\right) \varphi
_{bcd}.
\end{equation*}%
This quantity does not vanish in general, so we can conclude that 
\begin{equation}
\nabla _{X}\varphi \in \Lambda _{7}^{3}  \label{torsphi37}
\end{equation}%
and thus overall, 
\begin{equation}
\nabla \varphi \in W=\Lambda _{7}^{1}\otimes \Lambda _{7}^{3}.
\label{torsphiW}
\end{equation}%
Further classification of torsion classes depends on the decomposition of $W$
into components according to irreducible representations of $G_{2}$. Given (%
\ref{torsphiW}), we can write 
\begin{equation}
\nabla _{a}\varphi _{bcd}=T_{a}^{\ \ e}\psi _{ebcd}  \label{fulltorsion}
\end{equation}%
where $T_{ab}$ is the \emph{full torsion tensor}. This $2$-tensor full
defines $\nabla \varphi $ since pointwise, it has 49 components and the
space $W$ is also 49-dimensional (pointwise). In general we can split $%
T_{ab} $ as 
\begin{equation}
T=\tau _{1}g+\tau _{7}+\tau _{14}+\tau _{27}  \label{torsioncomps}
\end{equation}%
where $\tau _{1}$ is a function, and gives the $\mathbf{1}$ component of $T$%
, $\tau _{7}\in \Lambda _{7}^{2}$ and hence gives the $\mathbf{7}$
component, $\tau _{14}\in \Lambda _{14}^{2}$ gives the $\mathbf{14}$
component and $\tau _{27}$ is traceless symmetric, giving the $\mathbf{27}$
component. Note that the normalization of these components is different from 
\cite{karigiannis-2007}. Hence we can split $W$ as 
\begin{equation}
W=W_{1}\oplus W_{7}\oplus W_{14}\oplus W_{27}.  \label{Wsplit}
\end{equation}%
The 16 torsion classes arise as the subsets of $W$ which $\nabla \varphi $
belongs to. Moreover, as shown in \cite{karigiannis-2007}, the torsion
components $\tau _{i}$ relate directly to the expression for $d\varphi $ and 
$d\psi $. In fact, in our notation, 
\begin{subequations}%
\label{dpptors} 
\begin{eqnarray}
d\varphi &=&4\tau _{1}\psi +3\tau _{7}\wedge \varphi -\ast \tau _{27}
\label{dphitors} \\
d\psi &=&4\tau _{7}\wedge \psi -2\ast \tau _{14}.  \label{dpsitors}
\end{eqnarray}%
\end{subequations}%
Now suppose $d\varphi =d\psi =0$. Then this means that all four torsion
components vanish and hence $T=0$, and as a consequence $\nabla \varphi =0$.
The converse is trivially true, since $d$ and $d\ast $ can both be expressed
in terms of the covariant derivative. This result is due to Fern\'{a}ndez
and Gray \cite{FernandezGray}. If we add the fact that $Hol\left( g\right) $
is a subgroup of $G$ if and only if $X$ admits a torsion-free $G$ structure
from Proposition \ref{propholstruct}, then we get the following important
result.

\begin{theorem}[{Joyce, \protect\cite[Prop. 10.1.3]{Joycebook}}]
Let $X$ be a $7$-manifold with a $G_{2}$ structure defined by the $3$-form $%
\varphi $ and equipped with the associated Riemannian metric $g$. Then the
following are equivalent:\bigskip

\begin{enumerate}
\item The $G_{2}$-structure is torsion-free

\item $Hol\left( g\right) \subseteq G_{2}$ and $\varphi $ is the induced $3$%
-form

\item $\nabla \varphi =0$ on $X$ where $\nabla $ is the Levi-Civita
connection of $g$

\item $d\varphi =d\psi =0$ where $\psi =\ast \varphi $ with the Hodge star
defined by $g$
\end{enumerate}
\end{theorem}

Different torsion classes of the $G_{2}$ structure also restrict the
curvature of the manifold. Consider the curvature tensor $R_{abcd}.$ Then
for fixed $a$,$b$, we have 
\begin{equation*}
\left( R_{ab}\right) _{cd}\in \Lambda ^{2}\text{,}
\end{equation*}%
so we can decompose it as 
\begin{equation}
\left( R_{ab}\right) _{cd}=\left( \pi _{7}R_{ab}\right) _{cd}+\left( \pi
_{14}R_{ab}\right) _{cd}.  \label{riemcurvdecomp}
\end{equation}%
Following Karigiannis \cite{karigiannis-2007}, consider the operator $%
T_{\psi }$ (\ref{Tpsiop}) acting on $R_{abcd}$. Then we have 
\begin{eqnarray*}
g^{ad}T_{\psi }R_{abcd} &=&R_{abef}\psi _{\ \ \ cd}^{ef}g^{ad} \\
&=&-\left( R_{beaf}+R_{eabf}\right) \psi _{\ \ \ cd}^{ef}g^{ad} \\
&=&-R_{beaf}\psi _{\ \ c}^{e\ \ \ af}+R_{fbea}\psi _{\ \ \ \ c}^{eaf} \\
&=&-2g^{ad}T_{\psi }R_{abcd} \\
&=&0
\end{eqnarray*}%
where we have used the cyclic identity for $R_{abef}$. Hence, from (\ref%
{tpsil}) we get 
\begin{equation}
Ric_{bd}=3\left( \pi _{7}R_{ab}\right) _{cd}g^{ac}=\frac{3}{2}\left( \pi
_{14}R_{ab}\right) _{cd}g^{ac}  \label{riccidecom}
\end{equation}%
where $Ric_{bd}$ is the Ricci tensor. However, in general, by the
Ambrose-Singer holonomy theorem \cite{AmbroseSinger}, if $Hol\left( g\right)
\subseteq G$, then $R_{abcd}\in \mathrm{Sym}^{2}\left( \mathfrak{g}\right) $
where $\mathfrak{g}$ is the Lie algebra of $G$. Therefore, in the $G_{2}$
case, if the $G_{2}$ structure is torsion-free and hence $Hol\left( g\right)
\subseteq G_{2}$, then $R_{abcd}\in \mathrm{Sym}^{2}\left( \mathfrak{g}%
_{2}\right) $. This however implies that in (\ref{riemcurvdecomp}), the $\pi
_{7}$ component vanishes, and thus from (\ref{riccidecom}), we have the
following result:

\begin{theorem}[Bonan, \protect\cite{BonanG2}]
\label{propricciflat}Let $X$ be a Riemannian $7$-manifold with metric $g$.
If $Hol\left( g\right) \subseteq G_{2}$, then $X$ is Ricci-flat.
\end{theorem}

In fact, this result can also be derived without invoking the general
Ambrose-Singer theorem. In \cite{karigiannis-2007}, Karigiannis expressed
the $\Lambda _{7}^{2}$ component of the curvature tensor in terms of the
torsion tensor $T_{ab}$, so that when the torsion vanishes, the curvature
tensor is fully contained in $\Lambda _{14}^{2}$, thus directly confirming
the Ambrose-Singer theorem in the $G_{2}$ case. The original proof of
Theorem \ref{propricciflat} due to Bonan \cite{BonanG2} relied on the fact
that the Lie algebra structure of $\mathfrak{g}_{2}$ imposes strong
conditions on the Riemann tensor, and that these imply that the Ricci tensor
cannot be non-vanishing.

Given a compact manifold with a torsion-free $G_{2}$ structure, the
decompositions (\ref{formdecompose}) carry over to de Rham cohomology \cite%
{Joycebook}, so that we have 
\begin{subequations}%
\label{cohodecom} 
\begin{eqnarray}
H^{1}\left( X,\mathbb{R}\right) &=&H_{7}^{1} \\
H^{2}\left( X,\mathbb{R}\right) &=&H_{7}^{2}\oplus H_{14}^{2} \\
H^{3}\left( X,\mathbb{R}\right) &=&H_{1}^{3}\oplus H_{7}^{3}\oplus H_{27}^{3}
\\
H^{4}\left( X,\mathbb{R}\right) &=&H_{1}^{4}\oplus H_{7}^{4}\oplus H_{27}^{4}
\\
H^{5}\left( X,\mathbb{R}\right) &=&H_{7}^{5}\oplus H_{14}^{5} \\
H^{6}\left( X,\mathbb{R}\right) &=&H_{7}^{6}
\end{eqnarray}%
\end{subequations}%
Define the refined Betti numbers $b_{k}^{p}=\dim \left( H_{k}^{p}\right) $.
Clearly, $b_{1}^{3}=b_{1}^{4}=1$ and we also have $b^{1}=b_{7}^{k}$ for $%
k=1,...,6$. Moreover, it turns out that if $Hol\left( X,g\right) =G_{2}$
then $b^{1}=0$. Therefore, in this case the $H_{7}^{k}$ component vanishes
in (\ref{cohodecom}). It can be easily shown that on a Ricci-flat manifold,
any harmonic $1$-form must be parallel. However this happens if and only if $%
Hol\left( g\right) $ has an invariant $1$-form. However the only $G_{2}$%
-invariant forms are $\varphi $ and $\psi $. Therefore there are no
non-trivial harmonic $1$-forms when $Hol\left( g\right) =G_{2}$ and thus $%
b^{1}=0$.

An example of a construction of a manifold with a torsion-free $G_{2}$
structure is to consider $X=Y\times S^{1}$ where $Y$ is a Calabi-Yau $3$%
-fold. Define the metric and a $3$-form on $X$ as 
\begin{eqnarray}
g_{X} &=&d\theta ^{2}\times g_{Y}  \label{metCY} \\
\varphi &=&d\theta \wedge \omega +\func{Re}\Omega  \label{phiCY}
\end{eqnarray}%
where $\theta $ is the coordinate on $S^{1}$. This then defines a
torsion-free $G_{2}$ structure, with 
\begin{equation}
\ast \varphi =\frac{1}{2}\omega \wedge \omega -d\theta \wedge \func{Im}%
\Omega .  \label{psiCY}
\end{equation}%
However, the holonomy of $X$ in this case is $SU\left( 3\right) \subset
G_{2} $. From the K\"{u}nneth formula we get the following relations between
the refined Betti numbers of $X$ and the Hodge numbers of $Y$ 
\begin{eqnarray*}
b_{7}^{k} &=&1\ \ \ \text{for }k=1,...,6 \\
b_{14}^{k} &=&h^{1,1}-1\ \ \text{for }k=2,5 \\
b_{27}^{k} &=&h^{1,1}+2h^{2,1}\ \text{\ for }k=3,4.
\end{eqnarray*}

In \cite{JoyceG2} and \cite{Joycebook}, Joyce describes a possible
construction of a smooth manifold with holonomy equal to $G_{2}$ from a
Calabi-Yau manifold $Y$. So suppose $Y$ is a Calabi-Yau $3$-fold as above.
Then suppose $\sigma :Y\longrightarrow Y$ is an antiholomorphic isometric
involution on $Y$, that is, $\chi $ preserves the metric on $Y$ and
satisfies 
\begin{subequations}%
\begin{eqnarray}
\sigma ^{2} &=&1 \\
\sigma ^{\ast }\left( \omega \right) &=&-\omega \\
\sigma ^{\ast }\left( \Omega \right) &=&\bar{\Omega}.
\end{eqnarray}%
\end{subequations}%
Such an involution $\sigma $ is known as a \emph{real structure }on $Y$.
Define now a quotient given by 
\begin{equation}
Z=\left( Y\times S^{1}\right) /\hat{\sigma}  \label{barelydefine}
\end{equation}%
where $\hat{\sigma}$:$Y\times S^{1}\longrightarrow Y\times S^{1}$ is defined
by $\hat{\sigma}\left( y,\theta \right) =\left( \sigma \left( y\right)
,-\theta \right) $. The $3$-form $\varphi $ defined on $Y\times S^{1}$ by (%
\ref{phiCY}) is invariant under the action of $\hat{\sigma}$ and hence
provides $Z$ with a $G_{2}$ structure. Similarly, the dual $4$-form $\ast
\varphi $ given by (\ref{psiCY}) is also invariant. Generically, the action
of $\sigma $ on $Y$ will have a non-empty fixed point set $N$, which is in
fact a special Lagrangian submanifold on $Y$ \cite{Joycebook}. This gives
rise to orbifold singularities on $Z$. The singular set is two copies of $Z$%
. It is conjectured that it is possible to resolve each singular point using
an ALE $4$-manifold with holonomy $SU\left( 2\right) $ in order to obtain a
smooth manifold with holonomy $G_{2}$, however the precise details of the
resolution of these singularities are not known yet. We will therefore
consider only free-acting involutions, that is those without fixed points.

Manifolds defined by (\ref{barelydefine}) with a freely acting involution
were called \emph{barely }$G_{2}$ \emph{manifolds }by Harvey and Moore in 
\cite{Harvey:1999as}. The cohomology of barely $G_{2}$ manifolds is
expressed in terms of the cohomology of the underlying Calabi-Yau manifold $%
Y $: 
\begin{subequations}
\begin{eqnarray}
H^{2}\left( Z\right) &=&H^{2}\left( Y\right) ^{+} \\
H^{3}\left( Z\right) &=&H^{2}\left( Y\right) ^{-}\oplus H^{3}\left( Y\right)
^{+}
\end{eqnarray}%
\end{subequations}%
Here the superscripts $\pm $ refer to the $\pm $ eigenspaces of $\sigma
^{\ast }$. Thus $H^{2}\left( Y\right) ^{+}$ refers to two-forms on $Y$ which
are invariant under the action of involution $\sigma $ and correspondingly $%
H^{2}\left( Y\right) ^{-}$ refers to two-forms which are odd under $\sigma $%
. Wedging an odd two-form on $Y$ with $d\theta $ gives an invariant $3$-form
on $Y\times S^{1}$, and hence these forms, together with the invariant $3$%
-forms $H^{3}\left( Y\right) ^{+}$ on $Y$, give the three-forms on the
quotient space $Z$. Also note that $H^{1}\left( Z\right) $ vanishes, since
the $1$-form on $S^{1}$ is odd under $\hat{\sigma}$. Now, given a $3$-form
on $Y$, its real part will be invariant under $\sigma $, hence $H^{3}\left(
Y\right) ^{+}$ is essentially the real part of $H^{3}\left( Y\right) $.
Therefore the Betti numbers of $Z$ in terms of Hodge numbers of $Y$ are 
\begin{subequations}
\begin{eqnarray}
b^{1} &=&0 \\
b^{2} &=&h_{1,1}^{+} \\
b^{3} &=&h_{1,1}^{-}+h_{2,1}+1
\end{eqnarray}%
\end{subequations}%
A class of barely $G_{2}$ manifolds that are constructed from complete
intersection Calabi-Yau manifolds has recently been considered in \cite%
{GrigorianBarelyG2}, where the Betti numbers of all such manifolds have been
calculated explicitly.

Note that barely $G_{2}$ manifolds have holonomy $SU\left( 3\right) \ltimes 
\mathbb{Z}_{2}$ while the first Betti number still vanishes. This shows that
vanishing first Betti number is not a necessary and sufficient condition for 
$Hol\left( g\right) =G_{2}$. In fact, as shown by Joyce in \cite{JoyceG2}, $%
Hol\left( g\right) =G_{2}$ if and only if the fundamental group $\pi
_{1}\left( X\right) $ is finite.

Let us briefly describe Joyce's construction of compact torsion-free
manifolds with $Hol\left( g\right) =G_{2}$. Here we follow \cite{Joycebook}.
On $T^{7}$ we can define a flat $G_{2}$ structure $\left( \varphi
_{0},g_{0}\right) $, similarly as on $\mathbb{R}^{7}$. Now suppose that $%
\Gamma $ is a finite group acting on $T^{7}$ that preserves the $G_{2}$
structure. Then we can define the orbifold $T^{7}/\Gamma $. The key to
resolving the orbifold singularities is to consider appropriate \emph{Quasi
Asymptotically Locally Euclidean }(QALE)\emph{\ }$G_{2}$ manifolds. These
are $7$-manifolds with a torsion-free $G_{2}$ structure that is asymptotic
to the $G_{2}$ structure on $\mathbb{R}^{7}/G$ where $G$ is a finite
subgroup of $G_{2}$. The orbifold $T^{7}/\Gamma $ is then resolved to obtain
a smooth compact manifold. However on the resolution, the resulting $G_{2}$%
-structure is not necessarily torsion-free, so it is shown that it can be
deformed to a torsion-free $G_{2}$ structure $\left( \varphi ,g\right) $.
Further, the fundamental group is calculated, and if it is finite, then $%
Hol\left( g\right) =G_{2}$. Using this method, Joyce found 252 topologically
distinct $G_{2}$ holonomy manifolds with unique pairs of Betti numbers $%
\left( b^{2},b^{3}\right) $.

\section{Moduli space}

\setcounter{equation}{0}

\subsection{Deformations of $G_{2}$ structures}

\label{sectdeform}One of the interesting directions in the study of $G_{2}$
holonomy manifolds is the structure of the \emph{moduli space}. Essentially,
the idea is to consider the space of all torsion-free $G_{2}$ structures
modulo diffeomorphisms on a manifold with fixed topology. The moduli space
itself has an interesting geometry that may give further information about $%
G_{2}$ manifolds.

Currently, we can only say something about the very local structure of the $%
G_{2}$ moduli space. For this, we take a fixed $G_{2}$ structure and deform
it slightly. The space of these deformations is the local moduli space. To
study it, we thus need to understand the deformations of $G_{2}$ structures.
Although, we are mostly interested in deformations of torsion-free $G_{2}$
structures, many of the results are valid for any $G_{2}$ structures.

Our aim is to consider infinitesimal deformations of $\varphi $ of the form 
\begin{equation}
\varphi \longrightarrow \varphi +\varepsilon \chi  \label{g2deform}
\end{equation}%
for some $3$-form $\chi $. As we already know, the $G_{2}$ structure on $X$
and the corresponding metric $g$ are all determined by the invariant $3$%
-form $\varphi $. Hence, deformations of $\varphi $ will induce deformations
of the metric. These deformations of metric will then also affect the
deformation of $\psi =\ast \varphi $. Theoretically, \textquotedblleft
large\textquotedblright\ deformations could also be considered, and in fact,
as we shall see below in some cases closed expressions can be obtained for
large deformations. However in that case, it is difficult to determine the
resulting torsion class of the new $G_{2}$ structure \cite%
{karigiannis-2005-57}. In order for the deformed $\varphi $ to define a new $%
G_{2}$ structure, the new $\varphi $ must also be a positive form (as per
the definition of a $G_{2}$ structure). However it is known \cite{Joycebook}
that the bundle of positive $3$-forms on $X$ is an open subbundle of $%
\Lambda ^{3}T^{\ast }X$, so we can always find $\varepsilon $ small enough
in order for the deformed $\varphi $ to be positive.

Using the decomposition of $3$-forms (\ref{l3decom}), we can split $\chi $
into $\Lambda _{1}^{3},$ $\Lambda _{7}^{3}$ and $\Lambda _{27}^{3}$ parts,
and at first let us consider each one separately. As shown by Karigiannis in 
\cite{karigiannis-2005-57}, metric deformations can be made explicit when
the $3$-form deformations are either in $\Lambda _{1}^{3}$ or $\Lambda
_{7}^{3}$. Let us first review some of these results. First suppose 
\begin{equation}
\tilde{\varphi}=f\varphi  \label{phitildeconf}
\end{equation}%
We will also use the notation $\tilde{\psi}=\tilde{\ast}\tilde{\varphi}$
where $\tilde{\ast}$ is the Hodge star derived from the metric $\tilde{g}$
corresponding to $\tilde{\varphi}$. Then from (\ref{metriccomp1}) we get 
\begin{eqnarray}
\tilde{g}_{ab}\sqrt{\det \tilde{g}} &=&\frac{1}{144}\tilde{\varphi}_{amn}%
\tilde{\varphi}_{bpq}\tilde{\varphi}_{rst}\hat{\varepsilon}^{mnpqrst}  \notag
\\
&=&f^{3}g_{ab}\sqrt{\det g}  \label{p1g1}
\end{eqnarray}%
After taking the determinant on both sides, we obtain 
\begin{equation}
\det \tilde{g}=f^{\frac{14}{3}}\det g.  \label{p1detg1}
\end{equation}%
Substituting (\ref{p1detg1}) into (\ref{p1g1}), we finally get 
\begin{equation}
\tilde{g}_{ab}=f^{\frac{2}{3}}g_{ab}.  \label{p1g2}
\end{equation}%
and hence 
\begin{equation}
\tilde{\psi}=f^{\frac{4}{3}}\psi .  \label{p1sphi1}
\end{equation}%
So, a scaling of $\varphi $ gives a conformal transformation of the metric.
Hence deformations of $\varphi $ in the direction $\Lambda _{1}^{3}$ also
give infinitesimal conformal transformation. Suppose $f=1+\varepsilon a$,
then to third order in $\varepsilon $, we can write 
\begin{equation}
\tilde{\psi}=\left( \allowbreak 1+\frac{4}{3}a\varepsilon +\frac{2}{9}%
a^{2}\varepsilon ^{2}-\frac{4}{81}a^{3}\varepsilon ^{3}+O\left( \varepsilon
^{4}\right) \right) \psi \allowbreak .  \label{p1sphi2}
\end{equation}%
Given a torsion-free $G_{2}$ structure, $d\varphi =d\psi =0$, so if we want
the deformed structure to be also torsion-free, $f$ must be constant.

Now, suppose in general that $\tilde{\varphi}=\varphi +\varepsilon \chi $
for some $\chi \in \Lambda ^{3}$. Then using (\ref{metricdef}) for the
definition of the metric associated with $\tilde{\varphi}$, after some
manipulations, we get: 
\begin{eqnarray}
\widetilde{\left\langle u,v\right\rangle }\widetilde{\mathrm{vol}} &=&\frac{1%
}{6}\left( u\lrcorner \varphi \right) \wedge \left( v\lrcorner \varphi
\right) \wedge \varphi  \label{rhsdeform2} \\
&&+\frac{1}{2}\varepsilon \left[ \left( u\lrcorner \chi \right) \wedge \ast
\left( v\lrcorner \varphi \right) +\left( v\lrcorner \chi \right) \wedge
\ast \left( u\lrcorner \varphi \right) \right]  \notag \\
&&+\frac{1}{2}\varepsilon ^{2}\left( u\lrcorner \chi \right) \wedge \left(
v\lrcorner \chi \right) \wedge \varphi  \notag \\
&&+\frac{1}{6}\varepsilon ^{3}\left( u\lrcorner \chi \right) \wedge \left(
v\lrcorner \chi \right) \wedge \chi .  \notag
\end{eqnarray}%
Rewriting (\ref{rhsdeform2}) in local coordinates, we get%
\begin{equation}
\tilde{g}_{ab}\frac{\sqrt{\det \tilde{g}}}{\sqrt{\det g}}=g_{ab}+\frac{1}{2}%
\varepsilon \chi _{mn(a}\varphi _{b)}^{\ \ mn}+\frac{1}{8}\varepsilon
^{2}\chi _{amn}\chi _{bpq}\psi ^{mnpq}+\frac{1}{24}\varepsilon ^{3}\chi
_{amn}\chi _{bpq}\left( \ast \chi \right) ^{mnpq}  \label{rhsdeform3}
\end{equation}%
Now suppose the deformation is in the $\Lambda _{7}^{3}$ direction. This
implies that 
\begin{equation}
\chi =\omega \lrcorner \psi  \label{p7chi1}
\end{equation}%
for some vector field $\omega $. Look at the first order term in (\ref%
{rhsdeform3}). From (\ref{3formproj}) we see that this is essentially a
projection onto $\Lambda _{1}^{3}\oplus \Lambda _{27}^{3}$ - the traceless
part gives the $\Lambda _{27}^{3}$ component and the trace gives the $%
\Lambda _{1}^{3}$ component. Hence this term vanishes for $\chi \in \Lambda
_{7}^{3}$. For the third order term, it is more convenient to study at it in
(\ref{rhsdeform2}). By looking at 
\begin{equation*}
\omega \lrcorner \left( \left( u\lrcorner \omega \lrcorner \psi \right)
\wedge \left( v\lrcorner \omega \lrcorner \psi \right) \wedge \psi \right) =0
\end{equation*}%
we immediately see that the third order term vanishes. So now we are left
with%
\begin{eqnarray}
\tilde{g}_{ab}\sqrt{\det \tilde{g}} &=&\left( g_{ab}+\frac{1}{8}\varepsilon
^{2}\omega ^{c}\omega ^{d}\psi _{camn}\psi _{dbpq}\psi ^{mnpq}\right) \sqrt{%
\det g}  \notag \\
&=&\left( g_{ab}\left( 1+\varepsilon ^{2}\left\vert \omega \right\vert
^{2}\right) -\varepsilon ^{2}\omega _{a}\omega _{b}\right) \sqrt{\det g}
\label{p7gabtil}
\end{eqnarray}%
where we have used a contraction identity for $\psi $ twice. Taking the
determinant of (\ref{p7gabtil}) gives%
\begin{equation}
\sqrt{\det \tilde{g}}=\left( 1+\varepsilon ^{2}\left\vert \omega \right\vert
^{2}\right) ^{\frac{2}{3}}\sqrt{\det g}.  \label{p7detg}
\end{equation}%
Eventually we have the following result:

\begin{theorem}[Karigiannis, \protect\cite{karigiannis-2005-57}]
Given a deformation of a $G_{2}$ structure (\ref{g2deform}) with $\chi
=\omega \lrcorner \psi \in \Lambda _{7}^{3}$, then the new metric $\tilde{g}%
_{ab}$ is given by 
\begin{equation}
\tilde{g}_{ab}=\left( 1+\varepsilon ^{2}\left\vert \omega \right\vert
^{2}\right) ^{-\frac{2}{3}}\left( \left( g_{ab}\left( 1+\varepsilon
^{2}\left\vert \omega \right\vert ^{2}\right) -\varepsilon ^{2}\omega
_{a}\omega _{b}\right) \right)  \label{p7gab2}
\end{equation}%
and the deformed $4$-form $\tilde{\psi}$ is given by 
\begin{equation}
\tilde{\psi}=\left( 1+\varepsilon ^{2}\left\vert \omega \right\vert
^{2}\right) ^{-\frac{1}{3}}\left( \psi +\ast \varepsilon \left( \omega
\lrcorner \psi \right) +\varepsilon ^{2}\omega \lrcorner \ast \left( \omega
\lrcorner \varphi \right) \right) .  \label{p7starphi1}
\end{equation}
\end{theorem}

One of the key reasons why it is possible to get these closed form
expressions for modified $g$ and $\psi $ is because as shown by Karigiannis
in \cite{karigiannis-2005-57}, the determinant of (\ref{p7gabtil}) can be
calculated in a closed form. Notice that to first order in $\varepsilon $,
both $\sqrt{\det g}$ and $g_{ab}$ remain unchanged under this deformation.
Now let us examine the last term in (\ref{p7starphi1}) in more detail.
Firstly, we have 
\begin{equation*}
\omega \lrcorner \ast \left( \omega \lrcorner \varphi \right) =\ast \left(
\omega ^{\flat }\wedge \left( \omega \lrcorner \varphi \right) \right)
\end{equation*}%
and 
\begin{eqnarray}
\left( \omega ^{\flat }\wedge \left( \omega \lrcorner \varphi \right)
\right) _{mnp} &=&3\omega _{\lbrack m}\omega ^{a}\varphi _{\left\vert
a\right\vert np]}  \notag \\
&=&3\mathrm{i}_{\varphi }\left( \omega \circ \omega \right)  \label{omomphi1}
\end{eqnarray}%
where $\left( \omega \circ \omega \right) _{ab}=\omega _{a}\omega _{b}$.
Therefore, in (\ref{p7starphi1}), this term gives $\Lambda _{1}^{4}$ and $%
\Lambda _{27}^{4}$ components. So, can write (\ref{p7starphi1}) as 
\begin{equation}
\tilde{\psi}=\left( 1+\varepsilon ^{2}\left\vert \omega \right\vert
^{2}\right) ^{-\frac{1}{3}}\left( \left( 1+\frac{3}{7}\varepsilon
^{2}\left\vert \omega \right\vert ^{2}\right) \psi +\ast \varepsilon \left(
\omega \lrcorner \psi \right) +\varepsilon ^{2}\ast \mathrm{i}_{\varphi
}\left( \left( \omega \circ \omega \right) _{0}\right) \right) .
\label{p7starphi2}
\end{equation}%
Here $\left( \omega \circ \omega \right) _{0}$ denotes the traceless part of 
$\omega \circ \omega ,$ so that $\mathrm{i}_{\varphi }\left( \left( \omega
\circ \omega \right) _{0}\right) \in \Lambda _{27}^{3}$ and thus, in (\ref%
{p7starphi2}), the components in different representations are now
explicitly shown.

To first order, we thus have the deformations 
\begin{subequations}
\label{7def1stord}
\begin{eqnarray*}
\tilde{\varphi} &=&\varphi +\varepsilon \left( \omega \lrcorner \psi \right)
\\
\tilde{\psi} &=&\psi +\ast \varepsilon \left( \omega \lrcorner \psi \right) .
\end{eqnarray*}%
If originally $d\varphi =d\psi =0$, that is, the $G_{2}$ structure is
torsion-free, then for the deformed structure to be torsion-free to first
order we need 
\end{subequations}
\begin{equation*}
d\left( \omega \lrcorner \psi \right) =d\ast \left( \omega \lrcorner \psi
\right) =0.
\end{equation*}%
By expanding $d\left( \omega \lrcorner \psi \right) $ in terms of the
decomposition of $\Lambda ^{4}$, and setting each term individually to $0$,
we find that the symmetric part of $\nabla _{a}\omega _{b}$ and the $\Lambda
_{7}^{2}$ part of $d\omega ^{\flat }$ must vanish. Furthermore, by expanding 
$\ast d\ast \left( \omega \lrcorner \psi \right) $ in terms of the
decomposition of $\Lambda ^{2}$ we find that the $\Lambda _{14}^{2}$ part of 
$d\omega ^{\flat }$ must also vanish. Hence we get that $\nabla \omega =0$.
If $Hol\left( g\right) =G_{2}$, then we know that in this case $\omega =0$,
so there are no interesting small $\Lambda _{7}^{3}$ deformations of
manifolds with holonomy equal to $G_{2}$.

As we have seen above, in the cases when the deformations were in $\Lambda
_{1}^{3}$ or $\Lambda _{7}^{3}$ directions, there were some simplifications,
which make it possible to write down all results in a closed form. In the
case of deformations in $\Lambda _{27}^{3}$ the only known way to get
results for deformations of the metric and the $4$-form $\psi $ is to
consider the deformations order by order in $\varepsilon $. This analysis
has been carried out in \cite{GrigorianYau1}, and here we will review those
results. So suppose we have a deformation 
\begin{equation*}
\tilde{\varphi}=\varphi +\varepsilon \chi 
\end{equation*}%
where $\chi \in \Lambda _{27}^{3}$. Now let us set up some notation. Define 
\begin{eqnarray}
\tilde{s}_{ab} &=&\frac{1}{144}\frac{1}{\sqrt{\det g}}\tilde{\varphi}_{amn}%
\tilde{\varphi}_{bpq}\tilde{\varphi}_{rst}\hat{\varepsilon}^{mnpqrst}
\label{sabtilde} \\
&=&\tilde{g}_{ab}\sqrt{\frac{\det \tilde{g}}{\det g}}  \label{sabtilde2}
\end{eqnarray}%
From (\ref{metriccomp1}), the untilded $s_{ab}$ is then just equal to $g_{ab}
$. We can rewrite (\ref{sabtilde2}) as 
\begin{equation}
\tilde{g}_{ab}=\sqrt{\frac{\det g}{\det \tilde{g}}}\left( g_{ab}+\delta
s_{ab}\right)   \label{sabtilde3}
\end{equation}%
\qquad where $\delta g_{ab}$ is the deformation of the metric and $\delta
s_{ab}$ is the deformation of $s_{ab}$, which from (\ref{rhsdeform3}) is
given by%
\begin{equation}
\delta s_{ab}=\frac{1}{2}\varepsilon \chi _{mn(a}\varphi _{b)}^{\ \ mn}+%
\frac{1}{8}\varepsilon ^{2}\chi _{amn}\chi _{bpq}\psi ^{mnpq}+\frac{1}{24}%
\varepsilon ^{3}\chi _{amn}\chi _{bpq}\left( \ast \chi \right) ^{mnpq}.
\label{deltasab1}
\end{equation}%
Also introduce the following short-hand notation%
\begin{equation}
s_{k}=\func{Tr}\left( \left( \delta s\right) ^{k}\right)   \label{tkdef}
\end{equation}%
where the trace is taken using the original metric $g$. From (\ref{deltasab1}%
), note that since $\chi \in \Lambda _{27}^{3}$, when taking the trace the
first order term vanishes, and hence $s_{1}$ is at least second-order in $%
\varepsilon $. Clearly, for $k>1$, $s_{k}$ are at least of order $k$ in $%
\varepsilon $. Similarly as before, take the determinant of (\ref{sabtilde}):%
\begin{equation}
\left( \frac{\det \tilde{g}}{\det g}\right) ^{\frac{9}{2}}=\frac{\det \left(
g+\delta s\right) }{\det \left( g\right) }.  \label{detgl27}
\end{equation}%
Unlike in the case of $\Lambda _{7}^{3}$ deformations, we cannot compute $%
\det \left( g+\delta s\right) $ in closed form, so we have to calculate it
order by order in $\varepsilon $. From the standard expansion of $\det
\left( I+X\right) $, we find 
\begin{equation}
\frac{\det \left( g+\delta s\right) }{\det g}=1+s_{1}+\frac{1}{2}\left(
s_{1}^{2}-s_{2}\right) +\frac{1}{6}\left(
s_{1}^{3}-3s_{1}s_{2}+2s_{3}\right) +O\left( \varepsilon ^{4}\right) 
\label{detgtilde1}
\end{equation}%
However, as we noted above, $s_{1}$ is second-order in $\varepsilon $, so
this expression actually simplifies:%
\begin{equation}
\frac{\det \left( g+\delta s\right) }{\det g}=1+\left( s_{1}-\frac{1}{2}%
s_{2}\right) +\frac{1}{3}s_{3}+O\left( \varepsilon ^{4}\right) .
\label{detgtilde2}
\end{equation}%
Raising this to the power of $-\frac{1}{9}$, and expanding again to fourth
order in $\varepsilon $, we get 
\begin{equation}
\left( \frac{\det g}{\det \tilde{g}}\right) ^{\frac{1}{2}}=1+\left( \frac{1}{%
18}s_{2}-\frac{1}{9}s_{1}\right) -\frac{1}{27}s_{3}+O\left( \varepsilon
^{4}\right) .  \label{rdgtild5}
\end{equation}%
Using this and (\ref{sabtilde3}) we can immediately get the deformed metric,
but the expressions using the current form of $\delta s_{ab}$ are not very
useful. So far, the only property of $\Lambda _{27}^{3}$ that we have used
is that it is orthogonal to $\varphi $, thus in fact, up to this point
everything applies to $\Lambda _{7}^{3}$ as well. Now however, let $\chi $
be of the form 
\begin{equation}
\chi _{abc}=h_{[a}^{d}\varphi _{bc]d}  \label{chi27def}
\end{equation}%
where $h_{ab}$ is traceless and symmetric, so that $\chi \in \Lambda
_{27}^{3}$. Let us first introduce some further notation. Let $%
h_{1},h_{2},h_{3},h_{4}$ be traceless, symmetric matrices, and introduce the
following shorthand notation 
\begin{subequations}%
\begin{eqnarray}
\left( \varphi h_{1}h_{2}\varphi \right) _{mn} &=&\varphi
_{abm}h_{1}^{ad}h_{2}^{be}\varphi _{den}  \label{bphiab} \\
\varphi h_{1}h_{2}h_{3}\varphi  &=&\varphi
_{abc}h_{1}^{ad}h_{2}^{be}h_{3}^{cf}\varphi _{def}  \label{bphi} \\
\left( \psi h_{1}h_{2}h_{3}\psi \right) _{mn} &=&\psi _{abcm}\psi
_{defn}h_{1}^{ad}h_{2}^{be}h_{3}^{cf}  \label{bpsiab} \\
\psi h_{1}h_{2}h_{3}h_{4}\psi  &=&\psi _{abcm}\psi
_{defn}h_{1}^{ad}h_{2}^{be}h_{3}^{cf}h_{4}^{mn}  \label{bpsi}
\end{eqnarray}%
\end{subequations}%
It is clear that all of these quantities are symmetric in the $h_{i}$ and
moreover $\left( \varphi h_{1}h_{2}\varphi \right) _{mn}$ and $\left( \psi
h_{1}h_{2}h_{3}\psi \right) _{mn}$ are both symmetric in indices $m$ and $n$%
. Then, it can be shown that 
\begin{eqnarray*}
\chi _{(a\left\vert mn\right\vert }\varphi _{b)}^{\ \ mn} &=&\frac{4}{3}%
h_{ab} \\
\chi _{amn}\chi _{bpq}\psi ^{mnpq} &=&-\frac{4}{7}\left\vert \chi
\right\vert ^{2}g_{ab}+\frac{16}{9}\left( h^{2}\right) _{\{ab\}}-\frac{4}{9}%
\left( \varphi hh\varphi \right) _{\{ab\}} \\
\chi _{amn}\chi _{bpq}\ast \chi ^{mnpq} &=&\frac{32}{189}\func{Tr}\left(
h^{3}\right) g_{ab}-\frac{8}{9}\left( \varphi hh^{2}\varphi \right) _{\{ab\}}
\end{eqnarray*}%
where as before $\left\{ a\ b\right\} $ denotes the traceless symmetric
part. Using this and (\ref{deltasab1}), we can now express $\delta s_{ab}$
in terms of $h$:%
\begin{eqnarray}
\delta s_{ab} &=&\frac{2}{3}\varepsilon h_{ab}+g_{ab}\left( -\frac{1}{63}%
\varepsilon ^{2}\func{Tr}\left( h^{2}\right) +\frac{4}{567}\varepsilon ^{3}%
\func{Tr}\left( h^{3}\right) \right)   \label{deltasabfull} \\
&&+\varepsilon ^{2}\left( \frac{2}{9}\left( h^{2}\right) _{\{ab\}}-\frac{1}{%
18}\left( \varphi hh\varphi \right) _{\{ab\}}\right) -\frac{\varepsilon ^{3}%
}{27}\left( \varphi hh^{2}\varphi \right) _{\{ab\}}  \notag
\end{eqnarray}%
and hence 
\begin{subequations}%
\label{skfull} 
\begin{eqnarray}
s_{1} &=&\func{Tr}\left( \delta s\right) =-\frac{1}{9}\varepsilon ^{2}\func{%
Tr}\left( h^{2}\right) +\frac{4}{81}\varepsilon ^{3}\func{Tr}\left(
h^{3}\right)   \label{s1full} \\
s_{2} &=&\func{Tr}\left( \delta s^{2}\right) =\frac{4}{9}\varepsilon ^{2}%
\func{Tr}\left( h^{2}\right) +\varepsilon ^{3}\left( \frac{8}{27}\func{Tr}%
\left( h^{3}\right) -\frac{2}{27}\left( \varphi hhh\varphi \right) \right) 
\label{s2full} \\
s_{3} &=&\func{Tr}\left( \delta s^{3}\right) =\frac{8}{27}\varepsilon ^{3}%
\func{Tr}\left( h^{3}\right)   \label{s4full}
\end{eqnarray}%
\end{subequations}%
Substituting these expressions into (\ref{rdgtild5}) and (\ref{sabtilde3}),
we can get the full expression for the deformed metric (up to third order in 
$\varepsilon $) and correspondingly the expression for the deformed $4$-form 
$\psi :$

\begin{theorem}[Grigorian and Yau, \protect\cite{GrigorianYau1}]
Given a deformation of a $G_{2}$ structure (\ref{g2deform}) with $\chi
_{abc}=h_{[a}^{d}\varphi _{bc]d}\in \Lambda _{27}^{3}$, then the new metric $%
\tilde{g}_{ab}$ is given to third order in $\varepsilon $ by 
\begin{eqnarray}
\tilde{g}_{ab} &=&\left( 1+\frac{1}{18}\varepsilon ^{2}\func{Tr}\left(
h^{2}\right) +\frac{1}{81}\varepsilon ^{3}\func{Tr}\left( h^{3}\right) -%
\frac{1}{243}\varepsilon ^{3}\left( \varphi hhh\varphi \right) \right)
g_{ab}+\frac{2}{3}\varepsilon h_{ab}  \label{metricdeform} \\
&&+\varepsilon ^{2}\left( \frac{2}{9}\left( h^{2}\right) _{\left( ab\right)
}-\frac{1}{18}\left( \varphi hh\varphi \right) _{ab}\right) +\frac{2}{81}%
\varepsilon ^{3}h_{ab}\func{Tr}\left( h^{2}\right) -\frac{\varepsilon ^{3}}{%
27}\left( \varphi hh^{2}\varphi \right) _{ab}+O\left( \varepsilon ^{4}\right)
\notag
\end{eqnarray}%
and correspondingly, the deformed $4$-form $\tilde{\psi}$ is given by 
\begin{eqnarray}
\tilde{\psi} &=&\psi -\varepsilon \ast \chi +\varepsilon ^{2}\left( -\frac{1%
}{189}\func{Tr}\left( h^{2}\right) \psi +\frac{1}{6}\ast \mathrm{i}_{\varphi
}\left( \left( \phi hh\phi \right) _{0}\right) \right)  \label{starphi27} \\
&&+\varepsilon ^{3}\left( -\frac{2}{1701}\left( \varphi hhh\varphi \right)
\psi -\frac{5}{108}\func{Tr}\left( h^{2}\right) \ast \chi +\frac{1}{18}\ast 
\mathrm{i}_{\varphi }\left( h_{0}^{3}\right) -\frac{1}{36}\ast \mathrm{i}%
_{\varphi }\left( \left( \psi hhh\psi \right) _{0}\right) +\frac{1}{324}%
\alpha \wedge \varphi \right)  \notag \\
&&+O\left( \varepsilon ^{4}\right)  \notag
\end{eqnarray}%
where $\left( \phi hh\phi \right) _{0}$, $h_{0}^{3}$ and $\left( \psi
hhh\psi \right) _{0}$ denote the traceless parts of $\left( \phi hh\phi
\right) _{ab}$, $\left( h^{3}\right) _{ab}$ and $\left( \psi hhh\psi \right)
_{ab},$ respectively, and 
\begin{equation}
\alpha _{a}=\psi _{amnp}\varphi _{rst}h^{mr}h^{ns}h^{pt}  \label{p27ua}
\end{equation}
\end{theorem}

In general if such a deformation is performed on a torsion-free $G_{2}$
structure, then it is not known what conditions must $h$ satisfy in order
for the torsion class to be preserved. If we restrict our analysis only to
first order deformations, then it is easier to see these conditions.

Suppose we have $d\varphi =d\psi =0$ and we apply a deformation (\ref%
{g2deform}) with $\chi =\mathrm{i}_{\varphi }\left( h\right) $ for traceless
and symmetric. Then to first order the conditions for $d\tilde{\varphi}=d%
\tilde{\psi}=0$ are 
\begin{equation*}
d\chi =d\ast \chi =0.
\end{equation*}%
Hence the deformation must be a form that is closed and co-closed. For a
compact manifold this is thus equivalent to $\chi $ being harmonic. We can
also find what this condition means in terms of $h$. By decomposing $d\chi $
into $\Lambda _{1}^{4}$, $\Lambda _{7}^{4}$ and $\Lambda _{27}^{4}$
components, we find that we must have 
\begin{subequations}
\begin{eqnarray}
\nabla _{r}h_{\ a}^{r} &=&0  \label{pi7h1} \\
\nabla _{m}h_{a(b}\varphi _{\ \ \ \ c)}^{ma} &=&0  \label{pi27h1}
\end{eqnarray}%
\end{subequations}%
Further, if we decompose $\ast d\ast \chi $ into $\Lambda _{7}^{2}$ and $%
\Lambda _{14}^{2}$ components, we again get(\ref{pi7h1}) and moreover get a
new constraint 
\begin{subequations}
\label{h1}
\begin{equation}
\nabla _{m}h_{a[b}\varphi _{\ \ \ \ c]}^{ma}=0  \label{p14h1}
\end{equation}%
Thus overall, for $h$ traceless and symmetric, $\chi =\mathrm{i}_{\varphi
}\left( h\right) $ being closed and co-closed is equivalent to 
\end{subequations}
\begin{equation*}
\nabla _{r}h_{\ a}^{r}=0\ \ \ \text{and \ \ }\nabla _{m}h_{ab}\varphi _{\ \
\ \ c}^{ma}=0.
\end{equation*}%
On a compact manifold $\chi $ being closed and co-closed is equivalent to $%
\chi $ being harmonic. It also turns out \cite{AcharyaGukov} that, if $\chi $
is defined as above, then 
\begin{equation*}
\Delta \chi =0\ \Longleftrightarrow \ \Delta _{L}h=0
\end{equation*}%
where $\Delta _{L}$ is the Lichnerowicz operator given by 
\begin{equation}
\Delta _{L}h_{ab}=\nabla ^{2}h_{ab}+2R_{acbd}h^{cd}.  \label{lichop}
\end{equation}%
Therefore to preserve the torsion-free $G_{2}$ structure, we have to limit
our attention to zero modes of the Lichnerowicz operator. Note that, to
linear order, traceless deformations of the metric which preserve the Ricci
tensor are also precisely the Lichnerowicz zero modes, and this is
consistent with (\ref{metricdeform}) where the linear term in the metric
deformation is proportional to $h$.

Let us compare what happens here to what happens on Calabi-Yau manifolds 
\cite{Candelas:1990pi}. In that case, deformations of the metric $\delta
g_{mn}$ split into deformations of mixed type $\delta g_{\mu \bar{\nu}}$ and
deformations of pure type $\delta g_{\mu \nu }$ and $\delta g_{\bar{\mu}\bar{%
\nu}}$. From the mixed type deformations we can define a real $\left(
1,1\right) $-form 
\begin{equation}
i\delta g_{\mu \bar{\nu}}dx^{\mu }\wedge dx^{\bar{\nu}}  \label{kahlerdef}
\end{equation}%
and given the holomorphic $3$-form $\Omega $, we can use the mixed type
deformation to define a real $\left( 2,1\right) $-form 
\begin{equation}
\Omega _{\kappa \lambda }^{\ \ \ \ \bar{\nu}}\delta g_{\bar{\mu}\bar{\nu}%
}dx^{k}\wedge dx^{\lambda }\wedge dx^{\bar{\mu}}\text{.}  \label{complexdef}
\end{equation}%
In order to preserve the Calabi-Yau structure, the metric deformation must
preserve the vanishing Ricci curvature, and hence $\delta g_{mn}$ must
satisfy the Lichnerowicz equation: 
\begin{equation*}
\Delta _{L}\delta g_{mn}=0
\end{equation*}%
However, the Lichnerowicz equation for $\delta g_{mn}$ becomes equivalent to
both the $\left( 1,1\right) $-form (\ref{kahlerdef}) and the $\left(
2,1\right) $-form (\ref{complexdef}) being harmonic. Note that the
definition (\ref{complexdef}) is very similar to $\chi
_{abc}=h_{[a}^{d}\varphi _{bc]d}$ in the $G_{2}$ case with $\varphi $
playing the role of $\Omega $ and $h$ the role of $\delta g_{\bar{\mu}\bar{%
\nu}}$.

\subsection{Geometry of the moduli space}

In the theory of Calabi-Yau moduli spaces, one of the key results is that
the local moduli space of complex structure deformations is isomorphic to an
open set in $H^{m-1,1}\left( X\right) $ where $X$ is a Calabi-Yau $m$-fold.
Moreover, as it has been shown by Tian and Todorov \cite{Tian, TodorovWP},
any infinitesimal deformation can be in fact lifted to a full deformation.
For the moduli spaces of $G_{2}$ manifolds however, we can only replicate
the results about the local moduli space. First let us define the moduli
space of torsion-free $G_{2}$ structures. Let $\mathcal{X}$ be the set of of
positive $3$-forms $\varphi \in \mathcal{P}^{3}X$ such that $d\varphi =d\ast
_{\varphi }\varphi =0$. Here we use $\ast _{\varphi }$to emphasize that the
Hodge star is defined using the $G_{2}$ holonomy metric that is defined by $%
\varphi $ itself. Then $\mathcal{X}$ gives the set of all $3$-forms that
correspond to oriented, torsion-free $G_{2}$ structures. However we do not
want to distinguish between $3$-forms that are related by a diffeomorphism.
Hence, let $\mathcal{D}$ be the group of all diffeomorphisms of $X$ isotopic
to the identity. This group then acts naturally on $3$-forms. The \emph{%
moduli space }of torsion-free $G_{2}$ structures is then defined as the
quotient $\mathcal{M}=\mathcal{X}/\mathcal{D}$. The key result by Joyce is
that $\mathcal{M}$ is locally diffeomorphic to an open set of $H^{3}\left( X,%
\mathbb{R}\right) $:

\begin{theorem}[Joyce, \protect\cite{JoyceG2}]
Define a map $\Xi :\mathcal{X\longrightarrow }H^{3}\left( X,\mathbb{R}%
\right) $ by $\Xi \left( \varphi \right) =\left[ \varphi \right] $. Then $%
\Xi $ is invariant under the action of $\mathcal{D}$ on $\mathcal{X}$.
Moreover, $\Xi $ induces a diffeomorphism between neighbourhoods of $\varphi 
\mathcal{D\in M}$ and $\left[ \varphi \right] \in H^{3}\left( X,\mathbb{R}%
\right) $.
\end{theorem}

Since the dimension of $H^{3}\left( X,\mathbb{R}\right) $ is $b^{3}\left(
X\right) $, this result implies that $\dim \mathcal{M}=b^{3}\left( X\right) $%
. The full proof of this result can be found either in \cite{JoyceG2} or 
\cite{Joycebook}. This result covers the basic local properties of the $%
G_{2} $ moduli space, but we do not yet know anything about the global
structure of $\mathcal{M}$. So anything we can say about the moduli space
only holds in a small neighbourhood.

Looking back at the study of Calabi-Yau moduli spaces, we know that the
complex structure moduli space admits a K\"{a}hler structure, and the K\"{a}%
hler structure moduli space admits a Hessian structure \cite{Candelas:1990pi}%
. It turns out that on the $G_{2}$ moduli space we can also define a Hessian
structure. First let us define the notion of a \emph{Hessian manifold} \cite%
{Shima}

\begin{definition}
Let $M$ be a smooth manifold and suppose $D$ is a flat, torsion-free
connection on $M$. A Riemannian metric $G$ on a flat manifold $\left(
M,D\right) $ is called \emph{Hessian} if $G$ can be locally expressed as 
\begin{equation}
G=D^{2}H  \label{hessgdef1}
\end{equation}%
that is, 
\begin{equation}
G_{ij}=\frac{\partial ^{2}H}{\partial x^{i}\partial x^{j}}  \label{hessgdef2}
\end{equation}%
where $\left\{ x^{1},...,x^{n}\right\} $ is an affine coordinate system with
respect to $D$. Then $H$ is called the Hessian potential.
\end{definition}

Note that this is the closest analogue to a K\"{a}hler structure that can be
defined on a real manifold. In fact, as shown by Shima \cite{Shima}, if we
define a complex structure on the manifold $TM$, then the straightforward
extension of $G$ onto $TM$ is K\"{a}hler if and only if $G$ is a Hessian
metric on $\left( M,D\right) $. Thus the complexification of a Hessian
manifold is K\"{a}hler.

In the case of the $G_{2}$ moduli space $\mathcal{M}$, we know that $%
\mathcal{M}$ is locally diffeomorphic to an open set in $H^{3}\left( X,%
\mathbb{R}\right) $. Suppose we choose a basis $\left[ \varphi _{0}\right]
,...,\left[ \varphi _{n}\right] $ on $H^{3}\left( X,\mathbb{R}\right) $
where $n=b^{3}\left( X\right) -1$. Taking the unique harmonic
representatives of the basis elements, we can expand $\varphi \in \mathcal{M}
$ as 
\begin{equation}
\varphi =\sum_{N=0}^{n}s^{N}\phi _{N}.  \label{phiansatz}
\end{equation}%
Since $H^{3}\left( X,\mathbb{R}\right) $ is a vector space, $s_{0},...,s_{n}$
give an affine coordinate system, which in turn defines a flat connection $%
D=d$ on $\mathcal{M}$. It is trivial to check that this connection is
well-defined \cite{karigiannis-2007a}.

In order to define a metric on $\mathcal{M}$, we have to choose a Hessian
potential function on $\mathcal{M}$. The only natural function on $\mathcal{M%
}$ is the volume function $V\left( \varphi \right) $ given by (\ref{phiwpsi}%
): 
\begin{equation*}
V\left( \varphi \right) =\frac{1}{7}\int_{X}\varphi \wedge \psi .
\end{equation*}%
Note that as before, $\psi =\ast _{\varphi }\varphi $ is itself a function
of $\varphi $. So we can consider $V$ or some function of $V$ as potential
candidates for a Hessian potential. Let us calculate the Hessian of $V$.
Note that under a scaling $s^{M}\longrightarrow \lambda s^{M}$, $\varphi $
scales as $\varphi \longrightarrow \lambda \varphi $ and from (\ref{p1sphi1}%
), $\ast \varphi $ scales as $\ast \varphi \longrightarrow \lambda ^{\frac{4%
}{3}}\ast \varphi $, and so $V$ scales as 
\begin{equation*}
V\longrightarrow \lambda ^{\frac{7}{3}}V.
\end{equation*}%
So $V$ is homogeneous of order $\frac{7}{3}$ in the $s^{M}$, and hence%
\begin{eqnarray*}
s^{M}\frac{\partial V}{\partial s^{M}} &=&\frac{7}{3}V \\
&=&\frac{1}{3}\int s^{M}\phi _{M}\wedge \ast \varphi
\end{eqnarray*}%
and thus, 
\begin{equation}
\frac{\partial V}{\partial s^{M}}=\frac{1}{3}\int \phi _{M}\wedge \ast
\varphi .  \label{dvdsm}
\end{equation}%
Using our results on deformations of $G_{2}$ structures from Section \ref%
{sectdeform}, we can deduce that 
\begin{equation}
\partial _{N}\left( \ast \varphi \right) =\frac{4}{3}\ast \pi _{1}\left(
\phi _{N}\right) +\ast \pi _{7}\left( \phi _{N}\right) -\ast \pi _{27}\left(
\phi _{N}\right) .  \label{sphi1der}
\end{equation}%
Hence differentiating (\ref{dvdsm}) again, we find that 
\begin{eqnarray}
\frac{\partial V}{\partial s^{M}\partial s^{N}} &=&\frac{4}{9}\int \pi
_{1}\left( \varphi _{M}\right) \wedge \ast \pi _{1}\left( \varphi
_{N}\right) +\frac{1}{3}\int \pi _{7}\left( \varphi _{M}\right) \wedge \ast
\pi _{7}\left( \varphi _{N}\right)  \label{HessV} \\
&&-\frac{1}{3}\int \pi _{27}\left( \varphi _{M}\right) \wedge \ast \pi
_{27}\left( \varphi _{N}\right)  \notag
\end{eqnarray}%
Note that in the case when $b^{1}\left( X\right) =0$ (which in particular is
true when $Hol\left( g\right) =G_{2}$), since $H_{7}^{3}=H^{1}$, the $%
H_{7}^{3}$ component of $H^{3}\left( X,\mathbb{R}\right) $ is empty.
Therefore, the second term in (\ref{HessV}) vanishes, and we find that the
signature of this metric is Lorentzian - $\left( 1,b_{3}-1\right) $. Up to a
constant factor, this definition of the moduli space metric has been been
used in mathematical literature - in particular by Hitchin in \cite%
{Hitchin:2000jd} and Karigiannis and Leung in \cite{karigiannis-2007a}.
However in physics literature, in particular by Beasley and Witten in \cite%
{WittenBeasley} and by Gutowski and Papadopoulos in \cite{Gutowski:2001fm},
the potential $K$ given by 
\begin{equation}
K=-3\log V  \label{kahlerpot}
\end{equation}%
has been used instead.

The motivation for using this modified potential is two-fold. Firstly, this
is more in line with the logarithmic K\"{a}hler potentials on Calabi-Yau
moduli spaces. Secondly, and perhaps most importantly is that the metric
that arises from this potential appears as the target space metric of the
effective theory in four dimensions when the action for the $11$-dimensional
supergravity is reduced to four dimensions on a $G_{2}$ manifold. We will
hence define the moduli space metric $G_{MN}$ as 
\begin{equation*}
G_{MN}=\frac{\partial ^{2}K}{\partial s^{M}\partial s^{N}}.
\end{equation*}%
Using the definition of $K$ and (\ref{HessV}), we get 
\begin{eqnarray}
\frac{\partial ^{2}K}{\partial s^{M}\partial s^{N}} &=&\frac{1}{V}\left(
\int \pi _{1}\left( \varphi _{M}\right) \wedge \ast \pi _{1}\left( \varphi
_{N}\right) -\int \pi _{7}\left( \varphi _{M}\right) \wedge \ast \pi
_{7}\left( \varphi _{N}\right) \right.   \label{d2kss} \\
&&\left. +\int \pi _{27}\left( \varphi _{M}\right) \wedge \ast \pi
_{27}\left( \varphi _{N}\right) \right)   \notag
\end{eqnarray}%
In this case, if $b^{1}\left( X\right) =0$, we get 
\begin{equation}
G_{MN}=\frac{1}{V}\int_{X}\phi _{M}\wedge \ast \phi _{N}.  \label{d2kzz}
\end{equation}%
This metric is then in fact Riemannian. In the physics setting, apart from
the $G_{2}$ $3$-form, there is another $3$-form $C$ and when the $11$%
-dimensional supergravity action is reduced to four dimensions, the
parameters of $\varphi $ and $C$ naturally combine to give a
complexification of the $G_{2}$ moduli space. The extension of the metric $%
G_{MN}$ to this complex space is then K\"{a}hler \cite%
{WittenBeasley,GrigorianYau1,Gutowski:2001fm}. However since the metric on
the complexified space does not depend on $C$, there is not much difference
in treating the moduli space as a complexified K\"{a}hler manifold or a real
Hessian manifold. Here we will treat $\mathcal{M}$ as a real Hessian
manifold.

Now that we have fixed a metric on $\mathcal{M}$, we can proceed to various
other geometrical quantities. For this we will need to use higher
derivatives of $\psi $. In what follows we will assume that $b^{1}\left(
X\right) =0$, so that there no harmonic forms in $H_{7}^{3}$. Let us
introduce local special coordinates on $\mathcal{M}$. Let $\phi
_{0}=a\varphi $ and $\phi _{\mu }\in \Lambda _{27}^{3}$ for $\mu
=1,...,b_{27}^{3}$, so that $s^{0}$ defines directions parallel to $\varphi $
and $s^{\mu }$ define directions in $H_{27}^{3}$. Then, from the
deformations of $\psi $ in Section \ref{sectdeform}, we can extract the
higher derivatives of $\psi $ in these directions: 
\begin{subequations}%
\label{sphiderivs} 
\begin{eqnarray}
\partial _{0}\partial _{0}\psi  &=&\frac{4}{9}a^{2}\psi \ \ \ \ \ \ \
\partial _{0}\partial _{0}\partial _{0}\psi =-\frac{8}{27}a^{3}\psi 
\label{sphider1} \\
\partial _{0}\partial _{\mu }\psi  &=&-\frac{1}{3}a\ast \phi _{\mu }\ \ \ \
\ \ \partial _{0}\partial _{0}\partial _{\mu }\psi =\frac{2}{9}a^{2}\ast
\phi _{\mu }\ \ \   \label{sphider2} \\
\partial _{\mu }\partial _{\nu }\psi  &=&-\frac{2}{189}\func{Tr}\left(
h_{\mu }h_{\nu }\right) \psi +\frac{1}{3}\ast \mathrm{i}_{\varphi }\left(
\left( \varphi h_{\mu }h_{\nu }\varphi \right) _{0}\right)   \label{sphider3}
\\
\partial _{0}\partial _{\mu }\partial _{\nu }\psi  &=&\frac{4}{567}a\func{Tr}%
\left( h_{\mu }h_{\nu }\right) \psi -\frac{2}{9}a\ast \mathrm{i}_{\varphi
}\left( \left( \varphi h_{\mu }h_{\nu }\varphi \right) _{0}\right) 
\label{sphider4} \\
\partial _{\mu }\partial _{\nu }\partial _{\kappa }\psi  &=&-\frac{5}{18}%
\func{Tr}\left( h_{\mu }h_{\nu }\right) \ast \phi _{\kappa }+\frac{1}{3}\ast 
\mathrm{i}_{\varphi }\left( (h_{\mu }h_{\nu }h_{\kappa })_{0}\right) 
\label{sphider5} \\
&&-\frac{1}{6}\ast \mathrm{i}_{\varphi }\left( \left( \psi h_{\mu }h_{\nu
}h_{\kappa }\psi \right) _{0}\right) -\frac{4}{567}\left( \varphi h_{\mu
}h_{\nu }h_{\kappa }\varphi \right) \psi   \notag
\end{eqnarray}%
\end{subequations}%
where $h_{\mu }$,$h_{\nu }$ and $h_{\kappa }$ are traceless symmetric
matrices corresponding to the $3$-forms $\phi _{\mu }$,$\varphi _{\nu }$ and 
$\phi _{\kappa }$, respectively. On a Hessian manifold, there is a natural
symmetric $3$-tensor given by the derivative of the metric, or equivalently
the third derivative of the Hessian potential. We will denote this tensor $%
A_{MNP}$. By analogy with similar quantities on Calabi-Yau moduli spaces,
this tensors is called the \emph{Yukawa coupling}. Using these expressions,
following \cite{GrigorianYau1} we can now write down all the components of $%
A_{MNR}$: 
\begin{subequations}%
\label{Amnkp1p27} 
\begin{eqnarray}
A_{000} &=&-14a^{3}  \label{a000} \\
A_{00\mu } &=&0  \label{a00m} \\
A_{0\mu \nu } &=&-\frac{2a}{V}\int \phi _{\mu }\wedge \ast \phi _{\nu
}=-2aG_{\mu \nu }  \label{a0mn} \\
A_{\mu \nu \rho } &=&-\frac{2}{27V}\int \left( \varphi h_{\mu }h_{\nu
}h_{\rho }\varphi \right) dV  \label{amnr}
\end{eqnarray}%
\end{subequations}%
The full Riemann curvature on a Hessian manifold is then defined by 
\begin{equation}
\mathcal{R}_{\ \ NPQ}^{M}=\frac{1}{4}\left( A_{\ \ QR}^{M}A_{\ \
NP}^{R}-A_{\ PR}^{M}A_{\ NQ}^{R}\right) .  \label{hessriemcurv}
\end{equation}%
Note that since the fourth derivative of $K$ is fully symmetric, the fourth
derivative terms vanish here. However, we can also define the \emph{Hessian
curvature }tensor by 
\begin{equation}
\mathcal{Q}_{KLMN}=\partial _{M}\partial _{N}\partial _{L}\partial
_{K}K-A_{KMR}A_{\ \ LN}^{R}.  \label{modcurv1}
\end{equation}%
This tensor is the equivalent of the K\"{a}hler curvature, and carries more
information than the actual Riemann tensor (\ref{hessriemcurv}). The Riemann
curvature tensor is obtained from $\mathcal{Q}$ by 
\begin{equation}
\mathcal{R}_{MNPQ}=\frac{1}{2}\left( \mathcal{Q}_{MNPQ}-\mathcal{Q}%
_{NMPQ}\right) .
\end{equation}%
From (\ref{sphiderivs}), we can calculate the fourth derivatives of $K$, and
hence get all the components of $\mathcal{Q}$:

\begin{theorem}[Grigorian and Yau, \protect\cite{GrigorianYau1}]
The components of the Hessian curvature tensor $\mathcal{\ Q}$ corresponding
to the metric (\ref{d2kzz}) on the local moduli space of torsion-free $G_{2}$
structures are given by: 
\begin{subequations}
\begin{eqnarray}
\mathcal{Q}_{0000} &=&14a^{4}  \label{R0000} \\
\mathcal{Q}_{000\mu } &=&0  \label{R000m} \\
\mathcal{Q}_{00\mu \nu } &=&2a^{2}G_{\mu \nu }  \label{R00mn} \\
\mathcal{Q}_{0\mu \nu \rho } &=&-A_{\mu \nu \rho }a  \label{R0mnr} \\
\mathcal{Q}_{\kappa \mu \nu \rho } &=&\frac{1}{3}\left( G_{\mu \nu
}G_{\kappa \rho }+G_{\mu \kappa }G_{\nu \rho }-\frac{5}{7}G_{\mu \rho
}G_{\kappa \nu }\right) -G^{\tau \sigma }A_{\mu \tau \rho }A_{\kappa \nu
\sigma }  \label{Rkmnr} \\
&&+\frac{1}{V}\int \left( -\frac{2}{27}\func{Tr}\left( h_{\kappa }h_{\mu
}h_{\nu }h_{\rho }\right) +\frac{1}{27}\left( \psi h_{\kappa }h_{\mu }h_{\nu
}h_{\rho }\psi \right) +\frac{5}{81}\func{Tr}\left( h_{(\kappa }h_{\mu
}\right) \func{Tr}\left( h_{\nu }h_{\rho )}\right) \right) \mathrm{vol} 
\notag
\end{eqnarray}%
\end{subequations}%
\end{theorem}

Let us look in more detail at the expression for $A_{\mu \nu \rho }$. If we
define $h_{\mu }^{a}=h_{\mu \ m}^{\ a}dx^{m}$, then we get 
\begin{equation}
A_{\mu \nu \rho }=-\frac{4}{9V}\int \varphi _{abc}h_{\mu }^{a}\wedge h_{\nu
}^{b}\wedge h_{\rho }^{c}\wedge \psi .  \label{amnryuk}
\end{equation}%
Expressions for the $G_{2}$ Yukawa coupling has been derived by different
authors - in particular by Lee and Leung, \cite{Lee:2002fa}, de Boer, Naqvi
and Shomer \cite{deBoer:2005pt}, and Karigiannis \cite{karigiannis-2007}.
Similarly, we can rewrite (\ref{Rkmnr}) as 
\begin{eqnarray}
\mathcal{Q}_{\kappa \mu \nu \rho } &=&\frac{1}{3}\left( G_{\mu \nu
}G_{\kappa \rho }+G_{\mu \kappa }G_{\nu \rho }-\frac{5}{7}G_{\mu \rho
}G_{\kappa \nu }\right) -G^{\tau \sigma }A_{\mu \tau \rho }A_{\kappa \nu
\sigma }  \label{hesscurv2} \\
&&+\frac{8}{9}\frac{1}{V}\int \psi _{abcd}h_{\kappa }^{a}\wedge h_{\mu
}^{b}\wedge h_{\nu }^{c}\wedge h_{\rho }^{d}\wedge \varphi   \notag \\
&&+\frac{1}{81}\frac{1}{V}\int \left( 5\func{Tr}\left( h_{(\kappa }h_{\mu
}\right) \func{Tr}\left( h_{\nu }h_{\rho )}\right) -6\func{Tr}\left(
h_{\kappa }h_{\mu }h_{\nu }h_{\rho }\right) \right) \mathrm{vol}  \notag
\end{eqnarray}%
As we have mentioned previously, by complexifying the $G_{2}$ moduli space,
it is possible to turn the Hessian structure into a K\"{a}hler structure.
Similarly, the Hessian curvature $\mathcal{Q}$ becomes K\"{a}hler curvature.
On Calabi-Yau manifolds, the complex structure moduli space is naturally a
complex manifold, and admits a K\"{a}hler structure, while the K\"{a}hler
structure moduli space is naturally a Hessian manifold, but can be
complexified to become K\"{a}hler itself. We compare the various quantities
on $G_{2}$ moduli space and on the Calabi-Yau complex structure moduli space
in Figure \ref{compare}. 
\begin{figure}
\begin{equation*}
\begin{tabular}{|l|l|l|}
\hline
\textbf{Quantity} & $G_{2}$\textbf{\ moduli in }$\Lambda _{27}^{3}$ & 
\textbf{Complex structure moduli} \\ \hline
Form & $\varphi $, $\psi $ & $\Omega $ \\ \hline
Deformation space & $H_{27}^{3}$ & $H^{\left( 2,1\right) }$ \\ \hline
Metric deformation & $\frac{2}{3}h_{\mu \nu }$ & $\delta g_{\bar{\mu}\bar{\nu%
}}$ \\ \hline
Form deformation & $\chi _{abc}=h_{[a}^{d}\varphi _{bc]d}$ & $\chi _{\alpha
\beta \bar{\gamma}}=-\frac{1}{2}\Omega _{\alpha \beta }^{\ \ \ \ \bar{\delta}%
}\delta g_{\bar{\gamma}\bar{\delta}}$ \\ \hline
K\"{a}hler potential & $K=-3\log \left( \int \varphi \wedge \psi \right) $ & 
$K=-\log \left( i\int \Omega \wedge \bar{\Omega}\right) $ \\ \hline
Moduli space metric & $G_{\mu \nu }=\frac{1}{V}\int \phi _{M}\wedge \ast
\phi _{N}$ & $G_{\mu \bar{\nu}}=-\frac{\int \chi _{\mu }\wedge \bar{\chi}_{%
\bar{\nu}}}{\int \Omega \wedge \bar{\Omega}}$ \\ \hline
Yukawa coupling & $A_{\mu \nu \rho }=-\frac{4}{9V}\int \varphi _{abc}h_{\mu
}^{a}\wedge h_{\nu }^{b}\wedge h_{\rho }^{c}\wedge \psi $ & $\kappa _{\mu
\nu \rho }=-\int \Omega _{\alpha \beta \gamma }\chi _{\mu }^{\alpha }\wedge
\chi _{\nu }^{\beta }\wedge \chi _{\rho }^{\gamma }\wedge \Omega $ \\ \hline
Curvature & $\mathcal{Q}_{\kappa \mu \nu \rho }$ as in (\ref{hesscurv2}) & $%
\begin{array}{c}
\mathcal{R}_{\mu \bar{\nu}\rho \bar{\sigma}}=G_{\mu \bar{\nu}}G_{\rho \bar{%
\sigma}}+G_{\mu \bar{\sigma}}G_{\rho \bar{\nu}} \\ 
-e^{2K_{C}}\kappa _{\mu \nu }^{\ \ \ \ \bar{\tau}}\kappa _{\bar{\nu}\bar{%
\sigma}\bar{\tau}}%
\end{array}%
$ \\ \hline
\end{tabular}%
\end{equation*}%
\caption{Comparison of $G\sb 2$ moduli space and Calabi-Yau complex structure moduli space}%
\label{compare}%
\end{figure}%
We can see that there are a number of similarities. This leads to a
speculation that perhaps the $G_{2}$ moduli space possesses more structures
than it is currently known. One of the key features of Calabi-Yau moduli
spaces is the \emph{special geometry}, that is, both have a line bundle
whose first Chern class coincides with the K\"{a}hler class \cite%
{Freed:1997dp,Strominger:1990pd}. From physics point of view, special
geometry relates to the effective theory having $\mathcal{N}=2$
supersymmetry. M-theory compactified on $G_{2}$ manifolds only gives $%
\mathcal{N}=1$ supersymmetry, so from this point of view it is perhaps
unlikely that the (complexified) $G_{2}$ moduli space would admit precisely
this structure. Moreover, it was shown by Alekseevsky and Cort\'{e}s in \cite%
{alekseevsky-2008} that a so-called \emph{special real }structure on a
Hessian manifold corresponds to special K\"{a}hler structure on the tangent
bundle. A special real manifold is a Hessian manifold on which the cubic
form $DG$ (with $D$ being the flat connection, and $G$ the Hessian metric)
is parallel with respect to $D$. In our terms, this would mean that the
derivative of the Yukawa coupling $A$ vanishes. This is a rather strong
condition which is not necessarily fulfilled in our case. So perhaps instead
there is some intermediate structure that could be defined on the $G_{2}$
moduli space or its complexification.

\section{Concluding remarks}

In this paper we have reviewed the developments in the study of $G_{2}$
moduli spaces. Currently only the local picture of the moduli space is
known, so in the future it is natural to try and obtain at least some
information on the global structure of the $G_{2}$ moduli space. On
Calabi-Yau manifolds, the extension to the global moduli space was
originally done by Tian and Todorov \cite{Tian,TodorovWP}. We have seen that
there are a number of similarities in the local structure of Calabi-Yau
moduli spaces and $G_{2}$ moduli spaces, so it is feasible that it could
also be possible to derive similar global properties of $G_{2}$ moduli
spaces. However torsion-free $G_{2}$ structures are very non-linear in some
aspects - in particular, the metric depends non-linearly on $\varphi $ and
hence the differential equation $\nabla \varphi =0$ for a torsion-free
structure is also non-linear. Therefore, it is not clear how to extend
infinitesimal deformations of a $G_{2}$ structure to large deformations,
apart from considering deformations order by order. However even such
expansions quickly get very complicated.

Another possible topic for study would be to further develop approaches to
mirror symmetry on $G_{2}$ holonomy manifolds \cite{Gukov:2002jv}. One
possible direction for further research is to look at $G_{2}$ manifolds in a
slightly different way. Suppose we have type $IIA$ superstrings on a
non-compact Calabi-Yau $3$-fold with a special Lagrangian submanifold which
is wrapped by a $D6$ brane which also fills $M_{4}$. Then, as explained in 
\cite{VafaKlemm:2001nx}, from the $M$-theory perspective this looks like a $%
S^{1}$ bundle over the Calabi-Yau which is degenerate over the special
Lagrangian submanifold, but this $7$-manifold is still a $G_{2}$ manifold.
The moduli space of this manifold will be then determined by the Calabi-Yau
moduli and the special Lagrangian moduli. This possibly could provide more
information about mirror symmetry on Calabi-Yau manifolds \cite%
{StromingerYau:1996it}.

One more direction is to look at $G_{2}$ manifolds with singularities. So
far in this work we have considered only smooth $G_{2}$ manifolds, however,
from a physical point of view, $G_{2}$ manifolds with singularities are even
more interesting, as they yield more realistic matter content \cite%
{AcharyaWitten}. Also, the moduli spaces which we studied are for manifolds
with fixed topology. By allowing topological transitions through
singularities \cite{CveticGibbons}, it may be possible to find some
relations between the different moduli spaces. Understanding these questions
would improve our grasp of both the geometry and physics of $G_{2}$ moduli
spaces and the interplay between them.

\bibliographystyle{jhep-a}
\bibliography{refs2}

\end{document}